\documentclass[preprint,review,12pt]{elsarticle}
\usepackage{amsmath}
\usepackage{graphicx}
\usepackage{hyperref}
\usepackage{nomencl}
\usepackage{comment}
\usepackage{caption}
\usepackage{subcaption}

\usepackage{xcolor}
\usepackage{multirow}
\usepackage{subcaption}
\usepackage{amsthm} 
\newtheorem{lemma}{Lemma}  
\usepackage{amsmath,amsfonts,amssymb}
\usepackage{hyperref}
\usepackage{geometry}
\usepackage{multicol}
\usepackage{graphics}
\usepackage{lipsum}
\usepackage{array}
\usepackage{algorithm}
\usepackage{algpseudocode}

\makenomenclature
\usepackage{etoolbox}
\renewcommand\nomgroup[1]{%
  \item[\bfseries
  \ifstrequal{#1}{A}{Indices}{%
  \ifstrequal{#1}{B}{Sets}{%
  \ifstrequal{#1}{C}{Constants}{%
  \ifstrequal{#1}{D}{Variables}{}}}}%
]}
\begin{document}
\begin{frontmatter}

\title{Electrification of Transportation: A Hybrid Benders/SDDP Algorithm for Optimal Charging Station Trading}

\author[CTU]{Farnaz Sohrabi}
\author[CP]{Mohammad Rohaninejad}
\author[CTU]{J\'{u}lius Bem\v{s}}
\author[CP]{Zden\v{e}k Hanz\'{a}lek}

\address[CTU]{Faculty of Electrical Engineering, Czech Technical University in Prague, Prague, Czech Republic}
\address[CP]{Industrial Informatics Department, Czech Institute of Informatics Robotics and Cybernetics, Czech Technical University in Prague, Prague, Czech Republic}

\begin{abstract}
This paper examines the electrification of transportation as a response to environmental challenges caused by fossil fuels, exploring the potential of battery electric vehicles and hydrogen fuel cell vehicles as alternative solutions. However, a significant barrier to their widespread adoption is the limited availability of charging infrastructure. Therefore, this study proposes the development of comprehensive charging stations capable of accommodating both battery and hydrogen vehicles to address this challenge.
The energy is purchased from the day-ahead and intraday auction-based electricity markets, where the electricity price is subject to uncertainty. Therefore, a two-stage stochastic programming model is formulated while the price scenarios are generated utilizing a k-means clustering algorithm. Given the complexity of the proposed model, an efficient solution approach is developed through the hybridization of the Benders decomposition algorithm and stochastic dual dynamic programming.
In the Benders master problem, day-ahead bidding variables are determined, whereas the Benders sub-problem addresses intraday bidding and charging station scheduling variables, employing stochastic dual dynamic programming to tackle its intractability.
Additionally, we transform the mixed integer linear program model of the second stage problem into a linear program, confirming its validity through KKT conditions.
Our model provides practical insights for making informed decisions in electricity markets based on sequential auctions.
While the bidding curves submitted to the day-ahead market remain unaffected by scenarios, those submitted to the intra-day market show dependence on fluctuations in day-ahead market prices.
\end{abstract}


\end{frontmatter}

\printnomenclature
\nomenclature[A]{$h, t$}{Index for hourly day-ahead and quarter-hourly intraday intervals.}
\nomenclature[A]{$\xi_{1}, \xi_{2}, \xi$}{Index for day-ahead, intraday, and total scenarios.}
\nomenclature[A]{$i$}{Index for electricity bidding price steps.}

\nomenclature[C]{$\pi_{\xi}$}{Probabilities associated with scenarios.}
\nomenclature[C]{$\lambda_{d}(h,\xi_1), \lambda_{i}(t,\xi_2)$}{Hourly day-ahead and quarter-hourly intraday market prices.}
\nomenclature[C]{$\eta_{b}, \eta_{h}, \eta_{e}$}{Efficiency of the battery, tank, and electrolyzer.}
\nomenclature[C]{$b_{l,min}, b_{l,max}$}{Lower and upper capacity for battery level.}
\nomenclature[C]{$b_{c,min}, b_{c,max}$}{Lower and upper charging capacity for the battery.}
\nomenclature[C]{$b_{d,min}, b_{d,max}$}{Lower and upper discharging capacity for the battery.}
\nomenclature[C]{$e_{p,min}, e_{p,max}$}{Lower and upper capacity for the electrolyzer.}
\nomenclature[C]{$y(i)$}{Electricity market bidding price steps.}
\nomenclature[C]{$h_{l,min}, h_{l,max}$}{Lower and upper capacity for tank level.}
\nomenclature[C]{$h_{c,min}, h_{c,max}$}{Lower and upper charging capacity for the tank.}
\nomenclature[C]{$h_{d,min}, h_{d,max}$}{Lower and upper discharging capacity for the tank.}
\nomenclature[C]{$\mathrm{H}$}{Higher heating value of hydrogen.}
\nomenclature[C]{$l_{e}(t), l_{h}(t)$}{Electricity and hydrogen loads of the charging station.}
\nomenclature[C]{$\lambda_{e}, \lambda_{h}$}{Electricity and hydrogen selling prices.}

\nomenclature[D]{$m_{d}(t,\xi_1), m_{i}(t,\xi)$}{Purchased power from the day-ahead and intraday markets.}
\nomenclature[D]{$v_{e}(t,\xi), v_{h}(t,\xi)$}{Direct dedicated electricity and hydrogen for loads.}
\nomenclature[D]{$b_{c}(t,\xi), b_{d}(t,\xi)$}{Charged and discharged power from the battery.}
\nomenclature[D]{$e_{p}(t,\xi)$}{Power consumed by the electrolyzer.}
\nomenclature[D]{$b_{l}(t,\xi), h_{l}(t,\xi)$}{Stored power and hydrogen amount in the battery and tank.}
\nomenclature[D]{$\rho_{d}(i,h), \rho_{i}(i,t, \xi_1)$}{Bidding volume in the day-ahead and intraday markets.}
\nomenclature[D]{$h_{c}(t,\xi), h_{d}(t,\xi)$}{Charged and discharged hydrogen amount from the tank.}
\nomenclature[D]{$u_{b}(t,\xi), u_{h}(t,\xi)$}{Operation mode for the battery and tank.}

\section{Introduction} \label{Introduction}
Currently, global warming is a major issue for humanity, largely driven by the emission of greenhouse gases, particularly carbon dioxide from burning fossil fuels in transportation, which accounted for 26\% of global greenhouse gas emissions in 2019 \citep{dastjerdi2023transient}.
Vehicle electrification is a promising method for reducing air pollution, involving the replacement of components and systems in vehicles that traditionally operate on fossil fuels with those powered by electricity such as battery electric vehicles and hydrogen fuel cell vehicles  \citep{elmasry2024electricity}.
Electric vehicles use batteries to store electrical energy, while hydrogen fuel cell cars use a tank to store hydrogen, which is then converted into electricity through a fuel cell. The hydrogen utilized in fuel cell cars is obtained through electrolysis, a process that splits water into hydrogen and oxygen \citep{sohrabi2021strategic}.

A hybrid charging station, which combines electric charging facilities and hydrogen refueling infrastructures, is gaining significant attention due to its versatility and feasibility, offering concurrent services to both electric and hydrogen vehicles \citep{cai2023hierarchical}. 
There are limited studies that address the integration of electric vehicle charging stations with hydrogen refueling stations. 
The study in \citep{fang2023optimal} examines the integration of photovoltaic systems and hydrogen energy storage for charging electric vehicles and refueling hydrogen-powered vehicles. It highlights benefits such as reduced reliance on fossil fuels and improved system stability, while also acknowledging challenges like high initial costs and the complexity of managing multiple energy sources.
The study detailed in \citep{sanchez2019methodology} presents a design method for hybrid charging stations that incorporate renewable energy sources and energy storage systems. It underscores advantages such as lower greenhouse gas emissions and increased energy efficiency, but also considers drawbacks like substantial initial investments and the complexity of integrating diverse energy systems.
The research in \citep{raj2022energy} investigates the power management and stability of a system integrating diesel generators, photovoltaic cells, fuel cells, and battery storage. It emphasizes benefits such as enhanced energy reliability and reduced emissions, while also addressing challenges like high costs and the need for advanced control systems.
The development of a stand-alone refueling station, designed and planned with a focus on risk management for both hydrogen and electric vehicles, is discussed in \citep{sriyakul2021risk}. A two-stage stochastic model for a multi-functional charging station connected to the grid, aiming to maximize profits, is developed in \citep{sohrabi2022optimal}. An energy management model for a solar parking lot with electric vehicle charging and hydrogen vehicle refueling, employing a Markov decision process, is introduced in \citep{guo2022energy}. An optimal energy management model for an all-in-one station, which includes charging and swapping for electric vehicles, as well as refueling for hydrogen vehicles, is proposed in \citep{cciccek2022optimal}. An on-grid hybrid station that combines hydrogen refueling and battery swapping using an innovative hybrid robust optimization method is introduced in \citep{xu2022robust}. Detailed numerical simulations are conducted in \citep{schroder2020optimization} to investigate a grid-connected system that enables both battery electric vehicle charging and hydrogen fuel cell vehicle refueling. Mehrjerdi \citep{mehrjerdi2019off} proposes an optimal off-grid charging station design for electric and hydrogen vehicles, primarily powered by solar energy with a diesel generator backup. An approach to optimize the operation and bidding strategy of a multi-product charging station through a two-stage stochastic program is proposed in \citep{sohrabi2023coordinated}. Wang et al. \citep{wang2020robust} introduce a model for robustly optimizing the design of a standalone charging station for hydrogen and electric vehicles, primarily powered by photovoltaic systems with a backup diesel generator. Various renewable energy sources are integrated in \citep{ampah2022electric} to power a hydrogen refueling and electricity charging station in an African country. An optimal operation plan for controllable assets within a micro-grid to manage electric vehicle charging stations and a hydrogen refueling station is outlined in \citep{massana2022multi}. A framework to facilitate optimal configuration planning of an electric-hydrogen micro-energy system with vehicle sharing stations is proposed in \citep{zeng2023optimal}.

As novel charging station solutions become increasingly popular, comprehending their incorporation into the larger framework of electricity markets becomes essential. Our attention is directed towards day-ahead and intraday markets utilizing a discrete auction mechanism.
Various short-term pricing mechanisms in European power systems involve decisions on bidding formats, pricing rules, energy product granularity, and market timing.
In the electricity market clearing process, the bidding format is crucial, structured with each bid containing a quantity and a price pair for each time interval \citep{silva2022short}.
In auction-based trading markets, participants use marginal uniform pricing, submitting orders with delivery date, direction, quantity, and price. The Power Exchange aggregates these orders, determining the market clearing price and quantity at the intersection of the generation and demand curves \citep{braun2018price}.
In European markets, both day-ahead and intraday auctions coexist, differing in energy product granularity. The day-ahead market trades in one-hour periods, while the intraday market deals with shorter intervals \citep{ocker2017german}.
Another distinguishing feature is market timing, with the day-ahead market accepting price-volume combinations around noon and announcing clearing prices later. In contrast, the intraday auction closes after the release of day-ahead clearing prices, with results published thereafter \citep{finnah2022integrated}.

In this study, we explore coordinated bidding in daily auction markets with the goal of maximizing profit for charging station owners. We focus on accurately modeling price uncertainties and address the challenge of scenario reduction by utilizing the k-means algorithm for clustering. The k-means algorithm, originally introduced by MacQueen in \citep{macqueen1967some}, is applied to annual real market prices for this purpose. 
K-means clustering is widely recognized as a standard technique supported by various implementations, known for its speed \citep{seljom2021stochastic}.
In practice, k-means clustering often uses Lloyd's algorithm \citep{lloyd1982least}. Initially, the number of clusters is determined, clusters are randomly generated, and scenarios are categorized based on their distance from these clusters. The clusters are then updated to minimize the distance between data points and cluster centers \citep{yaghoubi2021optimal}.

To address the uncertainty in electricity prices, stochastic programming is commonly employed in scheduling problems. Scenario-based methods are often used within stochastic optimization, necessitating a significant number of discrete scenarios. However, the exponential increase in the number of scenarios poses computational challenges. Therefore, we adopt an approach based on a combined Benders/SDDP decomposition algorithm, as proposed in \citep{rebennack2014generation}. The algorithm formulates day-ahead bidding decisions through the Benders master problem, while the sub-problem handles intraday bidding decisions and stochastic charging station scheduling problems utilizing stochastic dual dynamic programming (SDDP).
The SDDP is an efficient approach to tackle the resulting large-scale stochastic linear programming problem, combining sampling-based and scenario-based nested Benders decomposition methods.

Introduced by J.F. Benders, the Benders decomposition technique is tailored to handle problems involving complicated variables \citep{bnnobrs1962partitioning}. Originally designed for solving mixed-integer linear programming problems, this algorithm has evolved significantly through various extensions, enabling its adaptation to a wide range of applications, including multi-stage stochastic programming \citep{rahmaniani2017benders}.
The Benders partition algorithm employs a decomposition approach with two levels, known as the master and sub-problem, establishing iterative sequences between these levels to achieve the joint optimal solution \citep{alizadeh2015dynamic}.
The utilization of the Benders decomposition approach aims to enhance the tractability of the original model and reduce its computational workload \citep{rahmani2017strategic}.
Compared to alternatives such as metaheuristics, it offers advantages such as strong algebraic principles, proven convergence, and flexibility for adjusting the optimality gap, making it widely applicable \citep{pishvaee2014accelerated}.

The SDDP algorithm, introduced by Pereira and Pinto, is a methodology for solving multi-stage stochastic optimization problems \citep{pereira1991multi}. Introducing stochasticity into models significantly increases the difficulty of solving optimization problems. While it's feasible to directly solve problems with a small number of scenarios, as the number of scenarios increases, finding a solution necessitates the use of sampling techniques for stochastic programming, among which the SDDP method is prominent \citep{lei2024optimal}.
SDDP, a nested Benders decomposition algorithm, employs sampling techniques to overcome the curse-of-dimensionality inherent to stochastic problems \citep{rebennack2016combining}. It iterates between forward and backward passes: the forward pass samples a subset of scenarios and computes the optimal solution for each sample path independently. In the backward pass, starting from the final stage, the algorithm adds supporting hyperplanes to the approximate cost-to-go functions of the prior stage. 
This approach ensures manageable computational performance and guarantees convergence during the solution process, making it possible to efficiently solve optimization problems without excessive computational burden, resulting in stable and reliable outcomes \citep{ding2021multi}.
SDDP has found wide application in solving various real-world problems such as hydrothermal scheduling, and several important algorithms, including approximate dual dynamic programming, are based on its principles \citep{steeger2017dynamic}.

In summary, our study focuses on the decision-making challenges faced by the owner of hybrid electric and hydrogen charging station participating in day-ahead and intraday electricity auction markets. 
We employ the k-means algorithm for scenario generation by clustering electricity annual real market prices.
To address the volatility of market prices, we employ a two-stage stochastic optimization framework. Managing computational complexities, we utilize a hybrid Benders/SDDP decomposition algorithm within a two-stage stochastic programming model. Day-ahead bidding decisions are determined by the Benders master problem, while intraday bidding decisions and charging station scheduling issues are managed by the sub-problem, employing the stochastic dual dynamic programming approach.
Our methodology is adaptable to markets featuring two auctions, as demonstrated through an illustrative case study, providing valuable insights for decision-making in electricity markets.
The results demonstrate feasibility and suggest the need for various strategies to mitigate the effects of uncertainties. The primary contributions of our research in this paper include:

\begin{itemize}
\item Designing a decision-making framework tailored for a hybrid electricity/hydrogen charging station participating in day-ahead and intraday auction markets.
\item Creating a model based on scenarios and utilizing the k-means clustering algorithm to account for uncertainties in electricity prices at the charging station.
\item Implementing a Benders decomposition within a two-stage stochastic programming model to determine day-ahead bidding decisions in the master problem and to address intraday bidding decisions and charging station scheduling issues in the sub-problem.
\item Utilizing the stochastic dual dynamic programming approach to tackle the sub-problem associated with the second-stage optimization.
\end{itemize} 

The paper is organized as follows: Section \ref{Methodology} provides an overview of the issue, develops the problem formulation, explains the scenario generation method, and elaborates on the solution algorithm. Section \ref{Case study} presents the case study. Finally, Section \ref{Conclusion} summarizes the paper.

\section{Methodology} \label{Methodology}

\subsection{Problem statement} \label{Problem statement}
\begin{figure}[hbt!]
    \centering    \includegraphics[width=0.7\columnwidth]{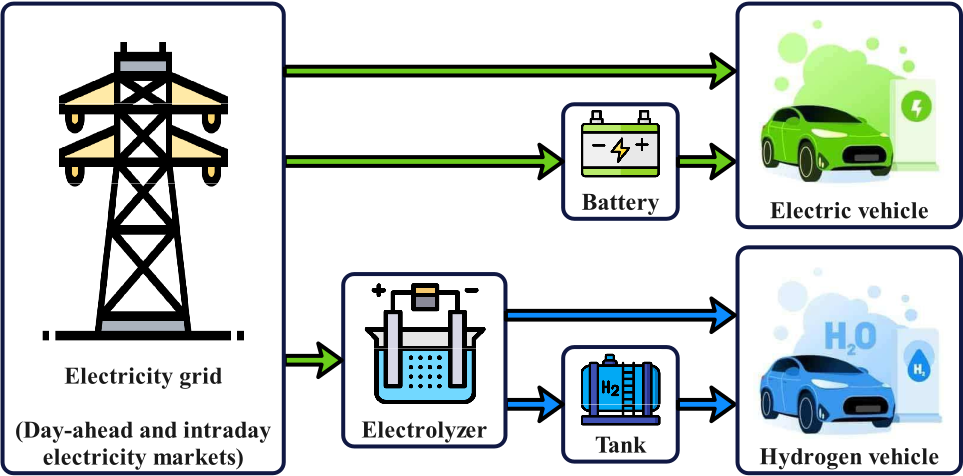}
    \caption{The design of our hybrid electric and hydrogen charging station.}
    \label{fig:station}
\end{figure}

In our study, we design a hybrid electric and hydrogen charging station, as illustrated in Fig. \ref{fig:station}. The purchased power by this grid-connected charging station is divided into three parts. The first component ensures the direct charging of electric vehicles, meeting their immediate power needs. The second component is dedicated to charging battery energy storage, contributing to the station's overall profitability. A significant aspect of this charging station is the third component, which contains an electrolyzer. This device converts electrical power into hydrogen, positioning hydrogen as a versatile energy carrier. The hydrogen generated in this process can be stored for future use or utilized for the efficient refueling of hydrogen vehicles, highlighting the station's multifaceted capabilities.

\begin{figure}[!hbt]
    \centering    \includegraphics[width=0.65\columnwidth]{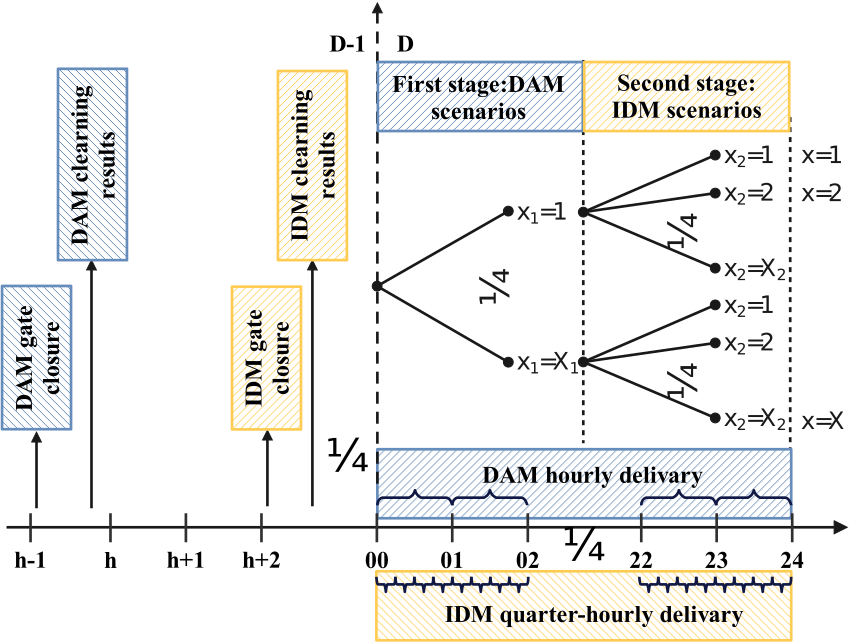}
    \caption{The sequence of events in the day-ahead market (DAM) and intraday market (IDM).}
    \label{fig:Markets}
\end{figure}

Our study explores the participation of charging station in both day-ahead and intraday auction markets, which are crucial components of electricity trading platforms. This involvement follows a structured sequence of actions, as depicted in Fig. \ref{fig:Markets}. The aim is to enhance the profitability of the charging station through the development of bidding decisions and the management of battery storage, electrolyzer, and hydrogen storage components. In the first stage, our focus is on devising precise day-ahead bidding curves. The progression to the second stage entails addressing a complex problem aimed at formulating intraday bidding curves and scheduling equipment operations.

The day-ahead auction market serves as the initial stage of decision-making. During this stage, charging station owners devise bidding curves for each hour of the following day. These curves are submitted before gate closure. It's important to note that all day-ahead bidding activities conclude after gate closure, signaling the start of the process to determine day-ahead market clearing prices.
Following this, our attention shifts to the intricacies of the intraday auction market. Here, the primary objective is to develop bidding curves for each intraday dispatch interval of the following day. The market time unit for European intraday auctions varies across countries, with intervals ranging from hourly to quarter-hourly. In our study, we concentrate on the quarter-hourly interval due to its complexity. By addressing the quarter-hourly scenario, we can subsequently handle half-hourly and hourly intervals. Consequently, we have a 96-step optimization horizon for each day. The submission of intraday bidding curves is permitted until the closure of the market, which occurs after the publication of day-ahead market clearing prices. Following the intraday gate closure, the clearing process commences, and intraday clearing prices are disclosed.

\begin{figure}[hbt!]
    \centering
    \includegraphics[width=0.45\columnwidth]{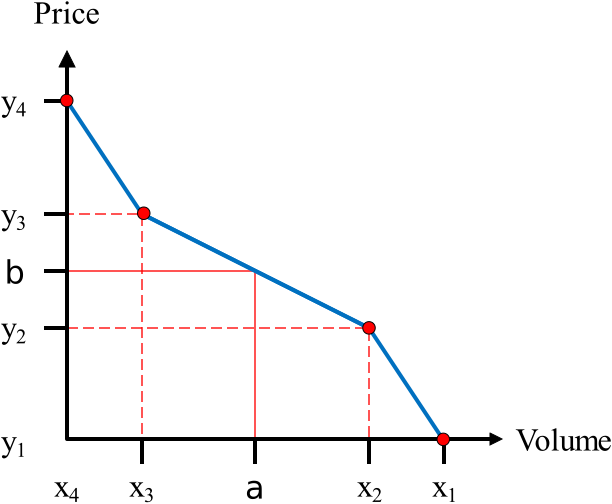}
    \caption{Bidding curve in the auction-based electricity market.}
    \label{fig:Bid curve}
\end{figure}

For each hour, bids are represented by points on a curve with a price and volume. The process involves choosing both, leading to a non-linear problem. To linearize it, prices $y$ are fixed, and volumes $x$ are decision variables. The curve is formed by linearly interpolating between price-volume pairs, creating a piece-wise linear curve \citep{fleten2005constructing}. Figure \ref{fig:Bid curve} illustrates a bidding curve with three segments and equations \eqref{eq.b1} and \eqref{eq.b2} represent the bidding curve with respect to prices and volumes, respectively.

\begin{flalign}\label{eq.b1}
    & \beta = y\left(i \right)+ \frac{y\left(i + 1\right)-y\left(i \right)}{x\left(i + 1\right)-x\left(i\right)} \left(\alpha -x\left(i\right) \right) \quad \text{if} \quad x\left(i + 1\right) \le \alpha \le x\left(i\right)  \notag \\
    & i \in \{1,...,I-1 \}&&
\end{flalign}
\begin{flalign}\label{eq.b2}
    & \alpha = \frac{\beta-y\left(i \right)}{y\left(i + 1\right)-y\left(i \right)} x\left(i + 1\right) + \frac{y\left(i + 1\right)-\beta}{y\left(i + 1\right)-y\left(i \right)} x\left(i\right) \quad \text{if} \quad y\left(i \right) \le \beta \le y\left(i + 1\right) , \notag \\
    & i \in \{1,...,I-1 \}&&
\end{flalign}

The probable difference in time frames between the two stages increases the overall complexity of this extensive stochastic optimization problem. The shift from day-ahead bidding decisions to intraday bidding decisions and operations requires an innovative approach capable of adapting to electricity market policies. Fundamentally, this constitutes the central challenge that our research aims to address.

\subsection{Mathematical modelling} \label{Mathematical modelling}
The sequential nature of the two-stage stochastic problem of optimal trading in electricity markets, where decisions unfold over time as more information becomes available through the stochastic process, is described in this section. 
In the first stage, decisions regarding day-ahead bidding are established. In the second stage, decisions regarding intraday bidding and operational actions are undertaken. Equation \eqref{eq.1} represents the objective function aiming to maximize the profit of the charging station. The cost of purchasing power from the day-ahead and intraday electricity markets is subtracted from the revenue generated by selling electricity and hydrogen to the vehicles to calculate the profit.
\begin{flalign}\label{eq.1}
    & \underset{\Xi}{\mathbf{maximize}}\Big[\sum_{t} 
    \big[l_{e}\left(t \right) \lambda_{e} + l_{h}\left(t \right) \lambda_{h}\big] -  \sum_{\substack{h, \xi_{1}, t \in h }} \pi_{\xi_1} \big [ m_{d}\left(t,\xi_1 \right) \lambda_{d}\left(h,\xi_1 \right)\big] - \notag \\
    &\sum_{\substack{t,\xi, \xi_{2} \in \xi}} \pi_{\xi} \big [ m_{i}\left(t,\xi \right) \lambda_{i}\left(t,\xi_2 \right)\big] \Big] \\
    &  \Xi = \big[m_{d}\left(t,\xi_1 \right), \rho_{d}\left(i,h \right) ,m_{i}\left(t,\xi \right),
    v_{e}\left(t,\xi\right),
    v_{h}\left(t,\xi\right),
    b_{l}\left(t,\xi\right), b_{c}\left(t,\xi\right), b_{d}\left(t,\xi\right), \notag \\
    &
    e_{p}\left(t,\xi\right), h_{l}\left(t,\xi\right), h_{c}\left(t,\xi\right), h_{d}\left(t,\xi\right),  \rho_{i}\left(i,t, \xi_1 \right),  u_{b}\left(t,\xi\right), u_{h}\left(t,\xi\right) \big] \notag&&
\end{flalign}

Equation \eqref{eq.2} enforces power balance across sequential markets, ensuring that the total power purchased from day-ahead and intraday markets equals the electricity designated for charging the battery and electric vehicles and running the electrolyzer.
\begin{flalign}\label{eq.2}
    m_{d}\left(t,\xi_1 \right) + m_{i}\left(t,\xi \right) =  v_{e}\left(t,\xi\right) + b_{c}\left(t,\xi\right) + e_{p}\left(t,\xi\right) \quad \forall t, \forall \xi, \forall \xi_{1} \in \xi &&
\end{flalign} 

Equation \eqref{eq.3} describes the power balance of the battery, ensuring that the battery's power level equals the sum of its previous level and the power purchased from markets for recharging the battery, minus the power discharged from the battery for electric vehicle refilling. Constraint \eqref{eq.4} bounds the inventory level of the battery, ensuring that the power level of the battery remains within its minimum and maximum capacities. Constraints \eqref{eq.5} and \eqref{eq.6} limit the charge and discharge capacity of the battery.
\begin{flalign}\label{eq.3}
b_{l}\left(t,\xi\right) = b_{l}\left(t-1,\xi\right) + \eta_{b}b_{c}\left(t,\xi\right) - b_{d}\left(t,\xi\right) / \eta_{b} \quad \forall t, \forall \xi &&
\end{flalign}
\begin{flalign}\label{eq.4}
    b_{l,min} \leq b_{l}\left(t,\xi\right) \leq b_{l,max}\quad \forall t, \forall \xi &&
\end{flalign}
\begin{flalign}\label{eq.5}
    & b_{c,min} u_{b}\left(t,\xi\right) \leq b_{c}\left(t,\xi\right) \leq b_{c,max} u_{b}\left(t,\xi\right) \quad \forall t, \forall \xi &&
\end{flalign}
\begin{flalign}\label{eq.6}
    & b_{d,min} \left( 1- u_{b}\left(t,\xi\right) \right) \leq b_{d}\left(t,\xi\right) \leq b_{d,max} \left( 1- u_{b}\left(t,\xi\right) \right) \quad \forall t, \forall \xi &&
\end{flalign}

Equation \eqref{eq.7h} ensures that the power purchased for hydrogen production, adjusted by the efficiency of the power-to-hydrogen conversion process and divided by the higher heating value of hydrogen, matches the allocated hydrogen for tank charging and refueling hydrogen vehicles. Constraint \eqref{eq.8h} imposes limits on the power consumption of the electrolyzer.
\begin{flalign}\label{eq.7h}
    & \eta_{e} e_{p}\left(t,\xi\right)/ \mathrm{H} = h_{c}\left(t,\xi\right)+ v_{h}\left(t,\xi\right)  \quad \forall t, \forall \xi &&
\end{flalign}
\begin{flalign}
    \label{eq.8h}
    & e_{p,min} \leq e_{p}\left(t,\xi\right) \leq e_{p,max} \quad \forall t, \forall \xi &&
\end{flalign}

Equation \eqref{eq.3h} outlines the hydrogen balance within the tank, ensuring that the hydrogen content in the tank matches the sum of its prior level and the hydrogen supplied by the electrolyzer for tank charging minus the hydrogen withdrawn from the tank for refueling hydrogen vehicles. Constraints \eqref{eq.4h} to \eqref{eq.6h} limit the inventory level, charging capacity, and discharging capacity of the hydrogen tank.
\begin{flalign}\label{eq.3h}
    & h_{l}\left(t,\xi\right) = h_{l}\left(t-1,\xi\right)+ \eta_{h} h_{c}\left(t,\xi\right) - h_{d}\left(t,\xi\right) / \eta_{h} \quad \forall t, \forall \xi &&
\end{flalign}
\begin{flalign}
    \label{eq.4h}
    & h_{l,min} \leq h_{l}\left(t,\xi\right) \leq h_{l,max} \quad \forall t, \forall \xi &&
\end{flalign}
\begin{flalign}
    \label{eq.5h}
    & h_{c,min} u_{h}\left(t,\xi\right) \leq h_{c}\left(t,\xi\right) \leq h_{c,max} u_{h}\left(t,\xi\right) \quad \forall t, \forall \xi &&
\end{flalign}
\begin{flalign}
    \label{eq.6h}
    & h_{d,min} \left( 1- u_{h}\left(t,\xi\right) \right) \leq h_{d}\left(t,\xi\right) \leq h_{d,max} \left( 1- u_{h}\left(t,\xi\right) \right) \quad \forall t, \forall \xi &&
\end{flalign}

Equation \eqref{eq.7} states that the power load for charging electric vehicles equals the sum of the power purchased from the markets and the power discharged from the batteries. Equation \eqref{eq.17h} ensures that the total hydrogen supplied from the electrolyzer and the hydrogen released from the tank for refueling hydrogen vehicles match the hydrogen demand.
\begin{flalign}\label{eq.7}
    v_{e}\left(t,\xi\right) + b_{d}\left(t,\xi\right) = l_{e}\left(t \right) \quad \forall t, \forall \xi &&
\end{flalign}
\begin{flalign}
    \label{eq.17h}
    & v_{h}\left(t,\xi\right) + h_{d}\left(t,\xi\right) = l_{h}\left(t \right) \quad \forall t, \forall \xi &&
\end{flalign}

Equation \eqref{eq.8} defines the hourly bidding curves in the day-ahead market as a piece-wise linear model described in \citep{fleten2005constructing}. Constraint \eqref{eq.9} enforces the market rules in the day-ahead market, stating that the hourly bidding curves should be non-increasing.
\begin{flalign}\label{eq.8}
    & \sum_{t \in h} m_{d}\left(t,\xi_1 \right) =\frac{\lambda_{d}\left(h,\xi_1 \right)-y\left(i \right)}{y\left(i+1 \right)-y\left(i \right)} \rho_{d}\left(i+1,h \right) + \frac{y\left(i+1 \right)-\lambda_{d}\left(h,\xi_1 \right)}{y\left(i+1 \right)-y\left(i \right)} \rho_{d}\left(i,h \right)\notag \\
    & if \quad y\left(i \right)\leq \lambda_{d}\left(h,\xi_1 \right) \leq y\left(i+1 \right) \quad  \forall h, \forall \xi_1, \forall i \in \{1,...,I-1 \} &&
\end{flalign}
\begin{flalign}\label{eq.9}
    \rho_{d}\left(i+1,h \right) \leq \rho_{d}\left(i,h \right) \quad \forall h, \forall i \in \{1,...,I-1 \} &&
\end{flalign}

Equation \eqref{eq.22} specifies the quarter-hourly bidding curves for the intraday market. Constraint \eqref{eq.24} ensures compliance with market rules, requiring that the quarter-hourly bidding curves be non-increasing in the intraday market. It's important to note that the bidding curves submitted to the intraday market depend on the outcomes observed in the day-ahead market scenarios. 
\begin{flalign}\label{eq.22}
    & m_{i}\left(t,\xi \right) =\frac{\lambda_{i}\left(t,\xi_2 \right)-y\left(i \right)}{y\left(i+1 \right)-y\left(i \right)} \rho_{i}\left(i+1,t, \xi_1 \right) + \frac{y\left(i+1 \right)-\lambda_{i}\left(t,\xi_2 \right)}{y\left(i+1 \right)-y\left(i \right)} \rho_{i}\left(i,t, \xi_1 \right) \notag\\
    & if \quad y\left(i \right)\leq \lambda_{i}\left(t,\xi_2 \right) \leq y\left(i+1 \right) \quad  \forall t, \forall \xi, \forall \left(\xi_1, \xi_2 \right) \in \xi, \forall i \in \{1,...,I-1 \}&&
\end{flalign}
\begin{flalign}\label{eq.24}
    \rho_{i}\left(i+1,t, \xi_1 \right) \leq \rho_{i}\left(i,t, \xi_1 \right) \quad \forall t, \forall \xi_1,\forall i \in \{1,...,I-1 \} &&
\end{flalign}
\begin{flalign*}
& \mathbf{variables:} \\
    & m_{d}\left(t,\xi_1 \right), 
    \rho_{d}\left(i,h \right), m_{i}\left(t,\xi \right),
    v_{e}\left(t,\xi\right),
    v_{h}\left(t,\xi\right),
    b_{l}\left(t,\xi\right), b_{c}\left(t,\xi\right), b_{d}\left(t,\xi\right),
    e_{p}\left(t,\xi\right), \notag\\
    & h_{l}\left(t,\xi\right), 
    h_{c}\left(t,\xi\right), h_{d}\left(t,\xi\right),  \rho_{i}\left(i,t, \xi_1 \right) \in \mathbb{R}\\
    & u_{b}\left(t,\xi\right), u_{h}\left(t,\xi\right) \in \{0, 1 \} &&
\end{flalign*}

\subsection{Scenario generation} \label{K-means method}
During the operation of the charging station in the short term, day-ahead and intraday electricity prices serve as input parameters. To enhance the optimization of decision variables, it is crucial to account for the stochastic nature of these inputs. The economic feasibility of coordinated bidding and operations significantly depends on accurately representing the uncertainties associated with electricity prices.
One commonly employed technique for scenario generation is the k-means algorithm, known for its simplicity, speed, and widespread application \citep{zeng2020scenario}. This algorithm partitions data into clusters, ensuring that objectives within each cluster share the highest similarity in terms of features and the least similarity with objectives in other clusters \citep{noorollahi2022scenario}.
However, before employing such methods, it is imperative to remove outliers from the data to ensure the accuracy of our analysis.

An efficient and rapid approach for outlier identification utilizes the $\mathcal{Z}$-score \citep{ashraf2017voltage}, as represented by Equation \eqref{eq.d1}, where $x_{i}$, $\mu$, and $\sigma$ denote the individual data point, mean, and standard deviation, respectively. Data points with $\lvert \mathcal{Z}_i \rvert$ less than or equal to the appropriate threshold are included in the new dataset, while those exceeding this threshold are omitted. This process enhances the robustness of data analysis by mitigating the influence of extreme values on statistical inferences.
\begin{flalign}\label{eq.d1}
& \mathcal{Z}i = \frac{ x_i - \mu }{\sigma} &&
\end{flalign}

After removing outliers, we apply the k-means algorithm. We define the number of clusters experimentally by considering a balance between complexity and accuracy. This algorithm aims to minimize the objective function, represented by Equation \eqref{eq.k1}, where $J$ denotes the objective function, $X_{i}$ represents the $i$-th observation, $C_{j}$ is the center of the $j$-th cluster, and $K$ indicates the number of clusters:
\begin{flalign}\label{eq.k1}
& J= \sum_{j=1}^{K} \sum_{i \in C_{j}} ||X_{i} - C_{j}||^{2}&&
\end{flalign}

The k-means algorithm proceeds through the following steps:
\begin{enumerate}
    \item Choose the number of clusters ($K$).
    \item Assign each item to the nearest cluster center.
    \item Update cluster centers and compute Eq. \eqref{eq.k1} for each cluster.
    \item Iterate through steps 2 and 3 until objects no longer switch clusters.
\end{enumerate}

By clustering data points based on their similarities and iteratively updating cluster centers, the k-means algorithm indirectly enables the derivation of scenario probabilities from the distribution of data points within each cluster. The probability of each cluster in a k-means clustering outcome is determined by dividing the number of data points assigned to that cluster by the total number of data points.

\subsection{Solution algorithm} \label{Solution algorithm}
To address a two-stage stochastic problem, a common approach is to solve a single large linear programming problem that encompasses all time steps and scenarios. However, the exponential growth of the scenario tree can quickly make the problem too complex to solve all at once. In such cases, decomposition techniques are employed.
To tackle this issue, we developed a Benders decomposition algorithm. In our approach, the day-ahead bidding decisions are determined in the Benders master problem (MP), while the intraday bidding and the charging station scheduling decisions are handled in the Benders sub-problem (SP). However, due to the high frequency of intraday market updates, occurring every 15 minutes in our case, the SP remains computationally intractable. To address this challenge, we integrated an SDDP algorithm within the Benders framework to efficiently solve the SP. The SDDP algorithm is employed to manage the large-scale stochastic linear programming problem inherent in the SP, allowing us to handle the frequent market updates effectively.
In the subsequent section, we provide a detailed explanation of the implementation of this hybrid Benders/SDDP algorithm. 

\subsubsection{Benders decomposition} \label{Benders decomposition}
Benders decomposition (BD) breaks down the main problem into two parts: the master problem and the sub-problem. These are then solved iteratively, with information exchanged between them to achieve an optimal solution. In each iteration of our proposed BD, the MP generates a solution and provides an upper bound (UB) on the objective function. The SP evaluates the MP solution and provides a lower bound (LB) based on the best feasible solution found. The UB is updated with cuts or constraints added to the MP, while the LB is updated with the objective value of the feasible solution from the SP. The iterative process continues, with the UB and LB progressively reducing the gap between them. The algorithm terminates when an optimal or near-optimal solution is determined.

\textit{A. Master problem}

The master problem in the first iteration of our proposed BD is as follows: in this model, the objective function \eqref{eq.1q} provides
the upper bound for the original model.
\begin{flalign}\label{eq.1q}
    & \alpha \left( \overline{m}_{d}\left(n, t,\xi_1 \right) \right)= \underset{\Xi_1}{\mathbf{maximize}}\Bigg[\sum_{t} 
    \left[l_{e}\left(t \right) \lambda_{e} + l_{h}\left(t \right) \lambda_{h}\right] -  \notag\\
    & \sum_{\substack{h, \xi_{1}, t \in h }} \pi_{\xi_1} \left [ m_{d}\left(t,\xi_1 \right) \lambda_{d}\left(h,\xi_1 \right)\right] -  \mathcal{Z} \Bigg] \\
    &  \Xi_1 = \big[m_{d}\left(t,\xi_1 \right), \rho_{d}\left(i,h \right), \mathcal{Z}  \big]\notag&&
\end{flalign}
\begin{flalign} \label{eq.32b}
    & \mathbf{subject~to:} \quad \eqref{eq.8}-\eqref{eq.9}&&
\end{flalign}
\begin{flalign*}
    & m_{d}\left(t,\xi_1 \right), \rho_{d}\left(i,h \right) \in \mathbb{R}, \mathcal{Z} \in \mathbb{R}^+ &&
\end{flalign*}

\textit{B. Sub-problem}

The sub-problem of our proposed BD is as follows: if the solution obtained from the MP satisfies the feasibility conditions of the sub-problem, this feasible solution establishes a lower bound for the original model. Equation \eqref{eq.2n} shows the power balance constraint in the sub-problem where the power purchased from the day-ahead market is fixed.
\begin{flalign}\label{eq.41nn}
    & \underset{\Xi_{2}}{\text{maximize}} \Big[-\sum_{\substack{t,\xi, \xi_{2} \in \xi}} \pi_{\xi} \big [ m_{i}\left(t,\xi \right) \lambda_{i}\left(t,\xi_2 \right)\big] \Big] \\
    &  \Xi_{2} = \big[m_{i}\left(t,\xi \right),
    v_{e}\left(t,\xi\right),
    v_{h}\left(t,\xi\right),
    b_{l}\left(t,\xi\right), b_{c}\left(t,\xi\right), b_{d}\left(t,\xi\right), 
    e_{p}\left(t,\xi\right), h_{l}\left(t,\xi\right), h_{c}\left(t,\xi\right), \notag \\
    &h_{d}\left(t,\xi\right),  \rho_{i}\left(i,t, \xi_1 \right),  u_{b}\left(t,\xi\right), u_{h}\left(t,\xi\right) \big] \notag&&
\end{flalign}
\begin{flalign} \label{eq.41nn1}
    & \mathbf{subject~to:} \quad \eqref{eq.3}-\eqref{eq.17h}, \eqref{eq.22}, \eqref{eq.24}&&
\end{flalign}
\begin{flalign}\label{eq.2n}
    \overline{m}_{d}\left(t,\xi_1 \right) + m_{i}\left(t,\xi \right) =  v_{e}\left(t,\xi\right) + b_{c}\left(t,\xi\right) + e_{p}\left(t,\xi\right) \quad \forall t, \forall \xi, \forall \xi_{1} \in \xi &&
\end{flalign} 
\begin{flalign*}
    & m_{i}\left(t,\xi \right),
    v_{e}\left(t,\xi\right),
    v_{h}\left(t,\xi\right),
    b_{l}\left(t,\xi\right), b_{c}\left(t,\xi\right), b_{d}\left(t,\xi\right),
    e_{p}\left(t,\xi\right),  h_{l}\left(t,\xi\right), 
    h_{c}\left(t,\xi\right), \notag\\
    & h_{d}\left(t,\xi\right),  \rho_{i}\left(i,t, \xi_1 \right) \in \mathbb{R}\\
    & u_{b}\left(t,\xi\right), u_{h}\left(t,\xi\right) \in \{0, 1 \} &&
\end{flalign*}

The binary variables, which involve the charging and discharging of the battery and hydrogen, categorize the SP problem as an MILP. To ensure the accuracy of dual variables, we apply a lemma 1 to convert the MILP formulation of SP into an equivalent LP.

\begin{figure}[hbt!]
    \centering
\includegraphics[width=0.7\columnwidth]{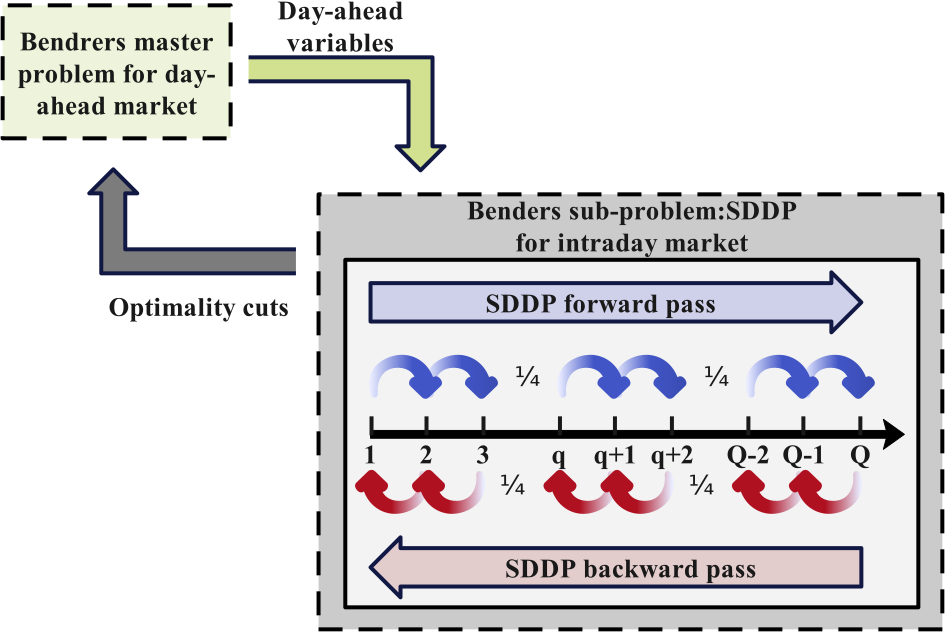}
    \caption{Forward and backward pass of SDDP.}
    \label{fig:SDDP}
\end{figure}

Since 96 quarter-hours make the SP computationally intractable, we divide it into $Q$ intervals, each containing a same number of quarter-hours, and employ SDDP, where $q \in \{1,...,Q\}$ represents the interval number. The SDDP framework decomposes the original problem into a series of interconnected secondary master problems and sub-problems, where each time interval $q$ is a sub-problem for the preceding interval $(q-1)$ and a master problem for the subsequent interval $(q+1)$. It involves an iterative process consisting of successive executions of the forward pass models $\{\mathcal{F}(q)|q=1,...,Q\}$ followed by the backward pass models $\{\mathcal{B}(q)|q=1,...,Q\}$, as illustrated in Fig. \ref{fig:SDDP}. This iterative approach aims to approximate profit functions efficiently.
The forward pass is executed $M$ times, each time randomly sampling a set $m$ $(m=1,...,M)$ of the price stochastic scenarios. The resulting variables are stored for the backward pass. 
During the backward pass, each $\mathcal{B}(q)$ is solved for the stored optimal state variables obtained in the forward pass and the scenario set $m$. In the backward pass, the algorithm calculates the expected cost of the $\mathcal{B}(q)$ and the dual variables, and the expected optimality cut is added to the $\mathcal{B}(q-1)$ consequently. 

At each interval \( q \), the secondary master problem \( \mathcal{F}(q) \) is formulated to optimize the immediate profit at the current interval and provide estimates of future profits based on decisions made at that interval. The corresponding secondary sub-problem \( \mathcal{F}(q+1) \) is formulated to optimize the profit at the next interval step, building upon decisions made in the current step. This iterative process, as illustrated in equation \eqref{eq.41}, ensures a coherent optimization strategy across the entire time horizon. The model \( \mathcal{F}(q) \) at iteration \( k+1 \) and \( \mathcal{B}(q) \) at iteration \( k \) in our proposed SDDP algorithm for the sampled scenario set \( m \) are defined as follows:
:
\begin{flalign}\label{eq.41}
    & \mathfrak{S}_{m}\left(\overline{b}_{l}\left(k, q \right), \overline{h}_{l}\left(k, q \right) \right) = \underset{\Xi_{2}}{\text{maximize}} \Bigg[ \left( -\sum_{t \in q,\xi' \in m, \xi_{2} \in \xi'} \left[ m_{i}(t,\xi') \lambda_{i}(t,\xi_2) \right] \right) /\mathcal{N}(m) - \notag\\
    &  \mathcal{S} \Bigg]&&
\end{flalign}
\begin{flalign} \label{eq.41b}
    & \mathbf{subject~to:} \quad \eqref{eq.4}-\eqref{eq.8h}, \eqref{eq.4h}-\eqref{eq.17h}, \eqref{eq.22}, \eqref{eq.24}&&
\end{flalign}
\begin{flalign}\label{eq.42}
    &  \mathcal{S} \geq - \mathfrak{S}_{m}\left(\overline{b}_{l}\left(k, q+1 \right), \overline{h}_{l}\left(k, q+1 \right) \right)  + \sum_{\substack{t \in q,\xi' \in m}}   \Big[ \big (\overline{b}_{l}\left(k, t,\xi' \right) - b_{l}\left(t,\xi' \right) \big) \tau_{b} \left(k, t,\xi' \right) + \notag \\ 
    & 
    \big (\overline{h}_{l}\left(k, t,\xi' \right) - h_{l}\left(t,\xi' \right) \big) \tau_{h} \left(k, t,\xi' \right)
    \Big] \quad  \forall k &&
\end{flalign}
\begin{flalign}\label{eq.32nc}
    & \overline{m}_{d}\left(t,\xi_1 \right) + m_{i}\left(t,\xi'\right) =  v_{e}\left(t,\xi'\right) + b_{c}\left(t,\xi'\right) + e_{p}\left(t,\xi'\right)  \quad  \forall t\in q, \forall \xi_{1} \in \xi', \notag \\ 
    & \forall \xi'\in m   &&
\end{flalign} 
\begin{flalign}\label{eq.3hn}
    & h_{l}\left(t,\xi'\right) = h_{l}\left(t-1,\xi'\right)+ \eta_{h} h_{c}\left(t,\xi'\right) - h_{d}\left(t,\xi'\right) / \eta_{h} \quad \forall t\in q, t \neq t_{in}\left(q\right),\notag \\ 
    &  \forall \xi' \in m &&
\end{flalign}
\begin{flalign}\label{eq.3n}
b_{l}\left(t,\xi'\right) = b_{l}\left(t-1,\xi'\right) + \eta_{b}b_{c}\left(t,\xi'\right) - b_{d}\left(t,\xi'\right) / \eta_{b} \quad \forall t\in q, t \neq t_{in}\left(q\right), \forall \xi' \in m &&
\end{flalign}
\begin{flalign}\label{eq.47}
    & \overline{b}_{l}\left(k, t_{f}\left(q-1\right),\xi' \right) = b_{l}\left(t_{in}\left(q\right),\xi'\right) - \eta_{b}b_{c}\left(t_{in}\left(q\right),\xi'\right) + b_{d}\left(t_{in}\left(q\right),\xi'\right) / \eta_{b} \notag \\ 
    &  : \tau_{b} \left(k, t,\xi' \right) \quad t = t_{in}\left(q\right), \forall \xi'\in m &&
\end{flalign}
\begin{flalign}\label{eq.48a}
    & \overline{h}_{l}\left(k, t_{f}\left(q-1\right),\xi' \right) = 
    h_{l}\left(t_{in}\left(q\right),\xi'\right)- \eta_{h} h_{c}\left(t_{in}\left(q\right),\xi'\right) + h_{d}\left(t_{in}\left(q\right),\xi'\right) / \eta_{h} \notag \\ 
    &  : \tau_{h} \left(k, t,\xi' \right) \quad t = t_{in}\left(q\right), \forall \xi'\in m &&
\end{flalign}
\begin{flalign*}
    & m_{i}\left(t,\xi' \right),
    v_{e}\left(t,\xi'\right),
    v_{h}\left(t,\xi'\right),
    b_{l}\left(t,\xi'\right), b_{c}\left(t,\xi'\right), b_{d}\left(t,\xi'\right),
    e_{p}\left(t,\xi'\right), h_{l}\left(t,\xi'\right), h_{c}\left(t,\xi'\right), \notag\\
    &h_{d}\left(t,\xi'\right),  \rho_{i}\left(i,t, \xi_1' \right) \in \mathbb{R}\\
    &  \mathcal{S} \in \mathbb{R}^+&&
\end{flalign*}

In this formula, \( \mathcal{N}(m) \) represents the number of scenarios within the chosen scenario sample \( m \), while \( t_{in}\left(q\right) \) and \( t_{f}\left(q-1\right) \) denote the first quarter in interval \( q \) and the last quarter in interval \( q-1 \), respectively.
Each secondary master problem is iteratively enriched with a series of optimality cuts, progressively refining the optimization process. The optimality cut described in equation \eqref{eq.42} restricts the potential future profits, with each inequality representing the impact of variations in the optimal value of the sub-problem concerning the state variables denoted as \( s_{t} \triangleq \{b_{l}(t,\xi'), h_{l}(t,\xi')\} \) within the secondary master problem. These state variables represent the state of the system at time \( t_{f}\left(q\right) \) and serve as crucial links between the master problem and the sub-problem. Subsequently, the optimal decisions for \( \mathcal{F}(q) \), derived from the state variable values obtained during the forward pass, are applied to \( \mathcal{F}(q+1) \). In each iteration, the effects of these decisions on the optimal value of the sub-problem \(\mathcal{F}(q+1)\) are captured by the corresponding dual variables \(\tau_{b}(k, t,\xi') \) and \( \tau_{h}(k, t,\xi')\).
Equations \eqref{eq.47} and \eqref{eq.48a} are utilized to compute the dual variables, which are then incorporated into all constraints within the secondary sub-problem \( \mathcal{B}(q) \) to generate the backward model \( \mathcal{B}(q-1) \). The primary objective is to guide the solution process towards state space regions that are likely to occur. Notably, the optimal state variables obtained in the forward path model \( \mathcal{F}(q) \) are stored for subsequent use in the backward pass to approximate \( \mathcal{B}(q) \).

Utilizing the dual variables \(\tau_{b}(k, t,\xi') \) and \( \tau_{h}(k, t,\xi')\), the optimality cuts \eqref{eq.42} are gradually constructed and added during each iteration of the backward pass, contributing to refining the approximation of profit functions.
The iterative process alternates between forward and backward passes until the profit function achieves a satisfactory level of precision. The algorithm stops when the difference between the upper and lower bounds of the optimal solution, calculated at each iteration, falls within a predefined tolerance value. These upper and lower bounds are determined by \eqref{eq.43} and \eqref{eq.44}, respectively. The convergence check is conducted using \eqref{eq.45} and \eqref{eq.46}, where $\sigma\left(k \right)$ denotes the standard deviation at iteration $k$ \citep{fatouros2017stochastic}.
\begin{flalign}\label{eq.43}
    & \overline{z}\left(k \right)=\mathcal{M}_{1} &&
\end{flalign}
\begin{flalign}\label{eq.44}
    & \underline{z}\left(k \right) = \Big[ -\sum_{\substack{t,\xi, \xi_{2} \in \xi}}  \big [ m_{i}\left(t,\xi \right) \lambda_{i}\left(t,\xi_2 \right)\big] \Big] / M &&
\end{flalign}
\begin{flalign}\label{eq.45}
    & \sigma\left(k \right) =\sqrt{\frac{1}{M-1}\sum_{m} \left(z\left(k \right) - \overline{z}\left(k \right)\right)^{2}  } &&
\end{flalign}
\begin{flalign}\label{eq.46}
    & \underline{z}\left(k \right) - 1.96 \sigma\left(k \right)\leq \overline{z}\left(k \right) \leq \underline{z}\left(k \right) + 1.96 \sigma\left(k \right)&&
\end{flalign}

\textit{C. Benders optimality cut}

Constraint \eqref{eq.32} illustrates the adjustments, referred to as optimality cuts, made using $\mathcal{Z}$ to approximate the sub-problem within the model. Equation \eqref{eq.32c} provides the dual values corresponding to the decisions made in the master problem. These dual variables, denoted as $\mu \left(m, n, t, \xi \right)$, are utilized in formulating the Benders optimality cuts for the next iteration. These cuts establish a connection between the master and sub-problems and are updated in each iteration to enhance the decisions of the master problem.
\begin{flalign}\label{eq.32}
    & \mathcal{Z} \geq - \mathfrak{S}_{m}\left(\overline{b}_{l}\left(k, 1,\xi' \right), \overline{h}_{l}\left(k, 1,\xi' \right) \right)   + \sum_{\substack{t,\xi, \xi_{1} \in \xi }} \Bigg[ \Big (\overline{m}_{d}\left(m, n, t,\xi_1 \right) - m_{d}\left(t,\xi_1 \right) \Big) \notag\\
    &  \mu \left(m, n, t,\xi \right) \Bigg] \quad \forall m, \forall n &&
\end{flalign}
\begin{flalign}\label{eq.32c}
    & \overline{m}_{d}\left(m, n, t,\xi_1 \right) =  v_{e}\left(t,\xi'\right) + b_{c}\left(t,\xi'\right) + e_{p}\left(t,\xi'\right) - m_{i}\left(t,\xi'\right) \quad :\mu \left(m, n, t,\xi' \right) \notag\\
    &\forall m, \forall n, \forall t,\forall \xi', \forall \xi_{1} \in \xi' &&
\end{flalign} 

\textit{D. Infeasibility state}

At times, the solution obtained from the Benders master problem fails to satisfy feasibility conditions due to the over-procuring of power from the day-ahead market. To address this issue, common approaches include incorporating slack variables into the sub-problem or adding feasibility cuts into the master problem. However, adding slack variables can cause numerical instability, and introducing feasibility cuts requires solving the sub-problem again to obtain the extreme ray, making the process time-consuming and complex. Therefore, these methods are not effective for our problem. Instead, we propose a more efficient approach by modifying equation \eqref{eq.7} to \eqref{eq.7n} to ensure sub-problem feasibility. This adaptation allows for potential over-procurement of power for direct electric vehicle charging while resulting in a feasible solution by acting as a penalty for over-purchasing in the objective function.
\begin{flalign}\label{eq.7n}
    v_{e}\left(t,\xi\right) + b_{d}\left(t,\xi\right) \geq l_{e}\left(t \right) \quad \forall t, \forall \xi &&
\end{flalign}


\textit{E. Outline of our algorithm}

\begin{figure}[hbt!]
    \centering
    \includegraphics[width=1\columnwidth]{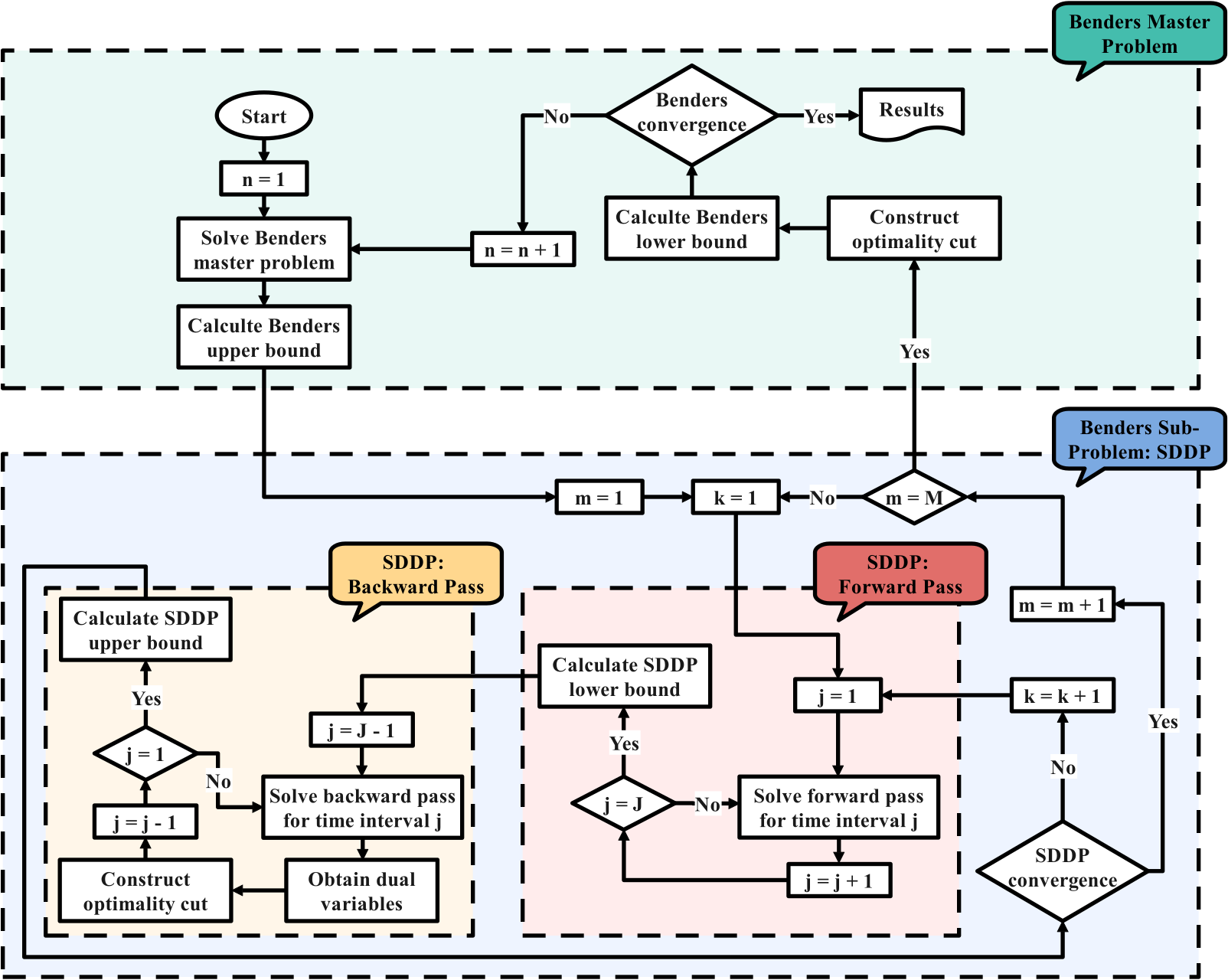}
    \caption{Solution algorithm for the optimal trading problem combining Benders decomposition with SDDP.}
    \label{fig:Solution}
\end{figure}

Our solution algorithm is depicted in Fig. \ref{fig:Solution} and Algorithm \ref{alg:hybrid_benders_sddp}. Here, $n$, $m$, $k$, and $j$ denote the number of Bender's iterations, SDDP samples, SDDP iterations, and SDDP time steps, respectively. Specifically, our algorithm involves \( n \) Bender's iterations. Within each Bender's iteration, we solve \( m \) SDDP samples, with each sample undergoing \( k \) iterations. During each SDDP iteration, both the forward and backward problems consist of \( j \) steps.

\begin{algorithm}
\caption{Hybrid Benders/SDDP Algorithm}
\label{alg:hybrid_benders_sddp}
\begin{algorithmic}[1]
\State Initialization
\While{$n < N$ \textbf{and} $|UB_{Benders} - LB_{Benders}| \leq \epsilon$}
    \State Solve Benders master problem; update $UB_{Benders}$
    \For{$m = 1$ to $M$}
        \While{$k < K$ \textbf{and} $|UB_{SDDP} - LB_{SDDP}| \leq \epsilon$}
            \For{$j = 1$ to $J$} 
                \State Solve forward pass for time interval $j$
            \EndFor
            \State Update $LB_{SDDP}$
            \For{$j = J$ to $1$}
                \State Solve backward pass for time interval $j$
                \State Construct SDDP optimality cut
            \EndFor
            \State Update $UB_{SDDP}$; $k \gets k + 1$
        \EndWhile
    \EndFor
    \State Add Benders cuts to MP; Update $LB_{Benders}$
    \State $n \gets n + 1$
\EndWhile
\end{algorithmic}
\end{algorithm}

\section{Results and discussion}\label{Case study}
The important details, like the prices for selling electricity and hydrogen, as well as the specifications of the battery storage, hydrogen tank, and electrolyzer, are listed in Table \ref{tab:data}.
\begin{table}[hbt!]
    \centering
    \caption{Required information of devices.}
    \label{tab:data}
    \begin{tabular}{l p{3cm} l l}
    \hline \hline
    Parameter & Value & Parameter & Value \\ 
    \hline
    $\lambda_{e}$ & $0.3 \; \text{(EUR/kWh)}$ & $\lambda_{h}$ & $12 \; \text{(EUR/kg)}$ \\
    $\eta_{b}$ & $0.85$ & $\eta_{h}$ & $0.9$ \\
    $b_{l,min}, b_{l,max}$ & $0, 60 \; \text{(kWh)}$  & $h_{l,min}, h_{l,max}$ & $0, 20 \; \text{(kg)}$ \\
    $b_{c,min}, b_{c,max}$& $0, 15 \; \text{(kW)}$ & $h_{c,min}, h_{c,min}$ & $0, 5 \; \text{(kg/h)}$ \\
    $b_{d,min}, b_{d,max}$& $0, 15 \; \text{(kW)}$ & $h_{d,min}, h_{d,max}$ & $0, 5 \; \text{(kg/h)}$ \\
    $e_{p,min}, e_{p,max}$& $0, 1000 \; \text{(kW)}$ & $\eta_{e}$& $0.8$ \\
    \hline \hline
    \end{tabular}
\end{table}

\begin{figure}[hbt!]
    \centering
\includegraphics[width=0.9\columnwidth]{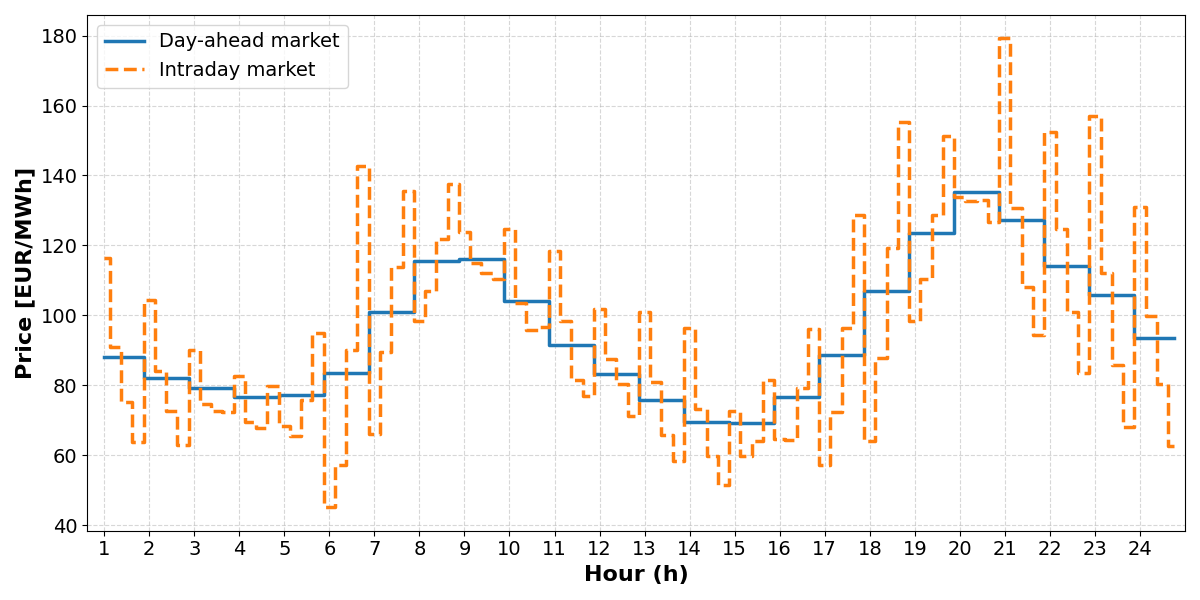}
    \caption{Day-ahead and intraday electricity market prices.}
    \label{fig:prices}
\end{figure}

The average hourly day-ahead and quarter-hourly intraday electricity prices are illustrated in Fig. \ref{fig:prices}.

\begin{figure}[hbt!]
    \centering
\includegraphics[width=0.9\columnwidth]{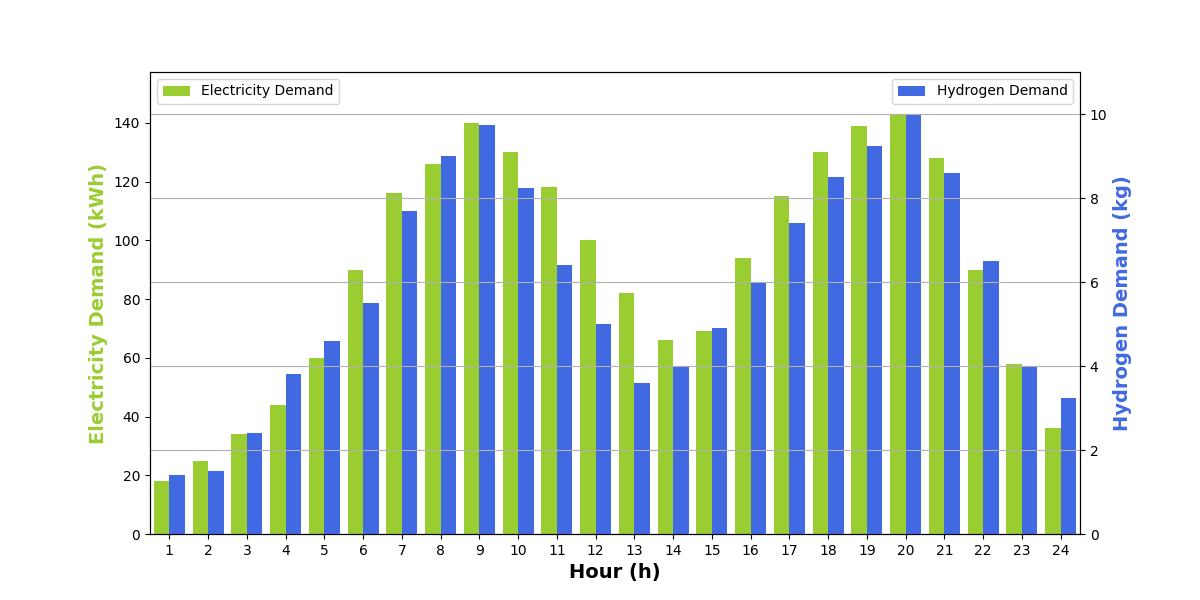}
    \caption{Electricity and hydrogen demand of the charging station.}
    \label{fig:demand}
\end{figure}

The electricity and hydrogen demand of vehicles are illustrated in Fig.~\ref{fig:demand} with green and blue bars, respectively. Although there are 96 quarters, we present the data on an hourly basis.

\begin{figure}[hbt!]
    \centering
    \subfloat[]{%
        \includegraphics[width=0.45\linewidth]{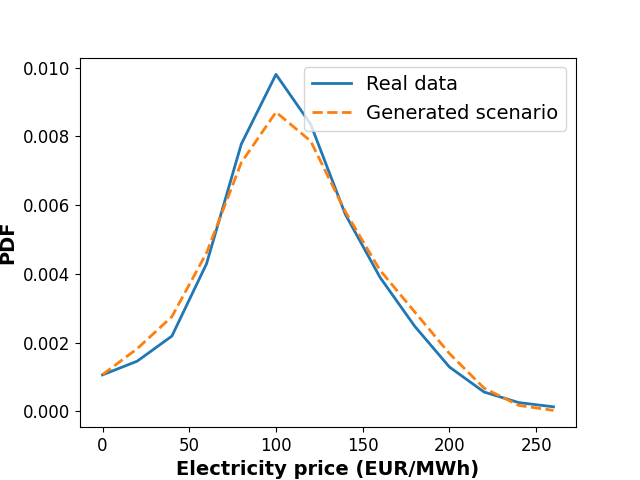}
        \label{fig:PDF_DM}
    }
    \quad 
    \subfloat[]{%
        \includegraphics[width=0.45\linewidth]{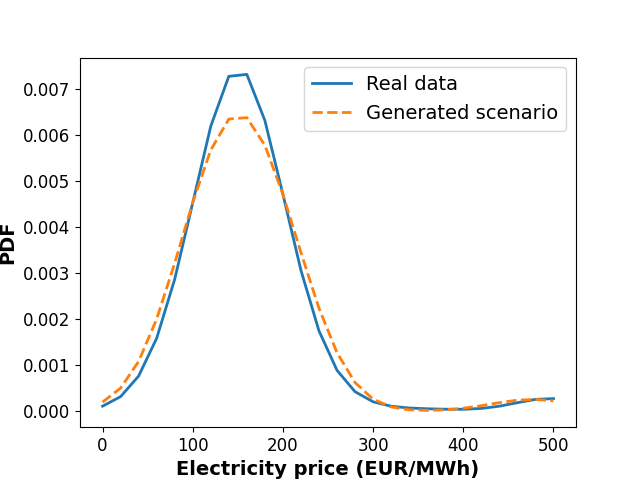}
        \label{fig:PDF_IM}
    }
    \caption{PDF of the original electricity price dataset versus generated scenarios by k-means. (a) Day-ahead market. (b) Intraday market.}
    \label{fig:PDF}
\end{figure}

Regarding the scheduling problem of the charging station, it is crucial to consider the uncertainty in electricity prices, as it can result in higher operational costs.
In order to assess the performance of our k-means approach for clustering electricity market price scenarios, we showcase the Probability Distribution Function (PDF) in Fig.~\ref{fig:PDF}. This kernel density function (KDE) plot incorporates scenarios generated by k-means, along with data sourced from the German market available on the ENTSO-E official website \citep{entsoe}. Our examination demonstrates a strong resemblance between the two PDFs, affirming the ability of our k-means model to effectively capture and replicate the intricate patterns inherent in real electricity market price dynamics.

Table \ref{tab:table1} illustrates the storage response to electricity price fluctuations over a two-hour period. When hourly day-ahead (DA) and quarter-hourly intraday (ID) electricity prices are low, the charging station strategically stores energy in both the battery and hydrogen tank, prioritizing direct provision to meet demand. Conversely, during intervals of high electricity prices, the station discharges stored energy from the battery and hydrogen tank, minimizing reliance on grid power. The charging station ensures not only the direct fulfillment of demand but also highlights its adaptability in managing storage units to optimize costs and align with market conditions, ultimately enhancing overall profitability.
\begin{table}[hbt!]
    \centering
    \caption{Charging Station Operation.}
    \label{tab:table1}
    \setlength{\tabcolsep}{3.5pt} 
    \begin{tabular}{l *{8}{c}}
        \hline \hline
        & \textbf{t1} & \textbf{t2} & \textbf{t3} & \textbf{t4} & \textbf{t5} & \textbf{t6} & \textbf{t7} & \textbf{t8} \\
        \hline 
        DA price (EUR/MWh) & 129 & 129 & 129 & 129 & 189 & 189 & 189 & 189 \\
        ID price (EUR/MWh) & 42 & 91 & 127 & 162 & 82 & 114 & 137 & 147 \\
        Electricity load (kWh) & 30 & 22 & 20 & 28 & 25 & 18 & 32 & 35 \\
        Hydrogen load (kg) & 5 & 10 & 12 & 8 & 15 & 6 & 14 & 9 \\
        Battery level (kWh) & 12.8 & 25.5 & 25.5 & 22.5 & 35.3 & 35.3 & 17.6 & 0 \\
        Hydrogen level (kg) & 4.5 & 9 & 9 & 6.6 & 11.1 & 11.1 & 5.6 & 0 \\
        \hline \hline
    \end{tabular}
\end{table}

\begin{figure}[hbt!]
    \centering
\includegraphics[width=0.6\columnwidth]{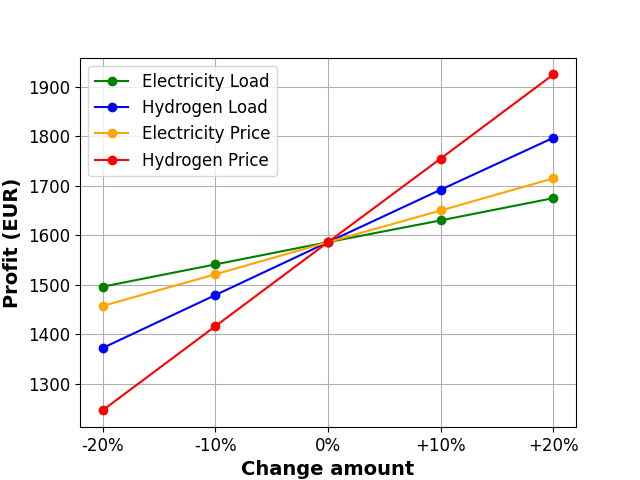}
    \caption{Sensitivity analysis.}
    \label{fig:sensitivity}
\end{figure}

The results of the profit sensitivity analysis reveal varying degrees of sensitivity among different parameters, including electricity load, hydrogen load, electricity price, and hydrogen price, as depicted in Fig. \ref{fig:sensitivity}. Notably, hydrogen price emerges as the most sensitive factor, demonstrating significant fluctuations in profit with its changes. Conversely, electricity load appears to be the least sensitive parameter, exhibiting relatively minimal variation in profit across the analyzed factors. The profit values for hydrogen load and electricity price fall somewhere in between, indicating a moderate sensitivity to changes. However, it's worth noting that profit is more sensitive to changes in hydrogen load compared to electricity price. These insights provide valuable information for decision-making, offering a better understanding of the impact of parameter variations on overall profit.

\begin{figure}[hbt!]
    \centering
    \subfloat[]{%
        \includegraphics[width=0.45\linewidth]{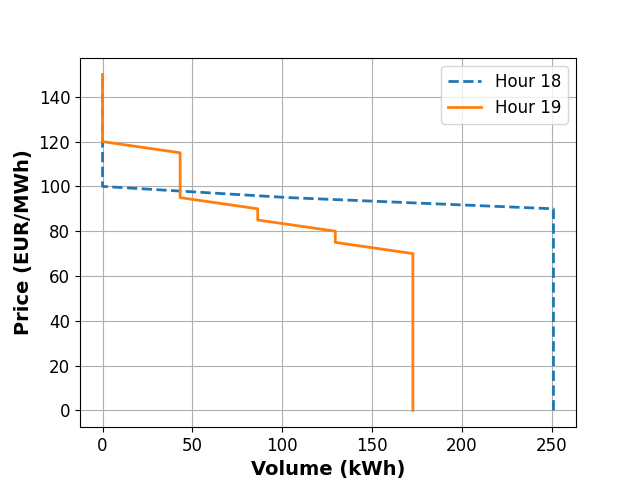}
        \label{fig:bid_DM}
    }
    \quad 
    \subfloat[]{%
        \includegraphics[width=0.45\linewidth]{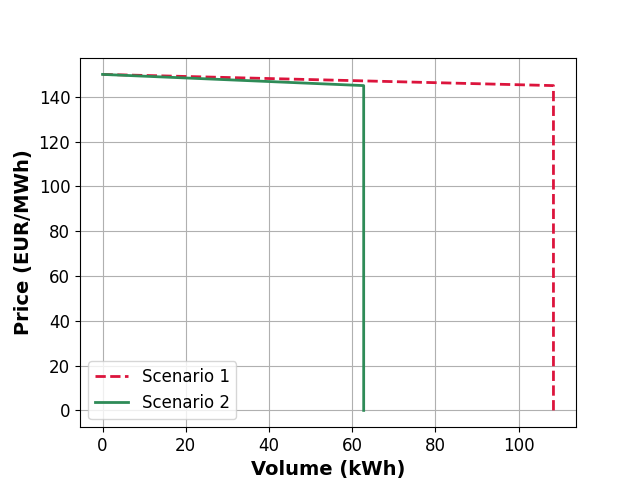}
        \label{fig:Bid_IM}
    }
    \caption{Bidding curves. (a) Day-ahead market. (b) Intraday market.}
    \label{fig:Bid}
\end{figure}

In Fig. \ref{fig:Bid} we try to show bidding curves in the day-ahead and intra-day markets to show that bidding in intra-day market not only depends on hours like day-ahead market but also depends on the day-ahead market biddding curves on each hour since submitted boidding curve to the intra-day market is the second satge variable and it is decided after the claerance of day-ahead market.
Fig. \ref{fig:bid_DM} displays day-ahead bidding curves for two selected hours within a 24-hour period analysis, aiming to highlight disparities in bidding behaviors across different time periods. Evidently, the power station procures higher volumes during hours with lower prices, reducing its purchases as prices rise. While day-ahead market bidding curves function as first stage variables and remain independent of scenarios, intra-day market bidding curves act as second stage variables, depending on scenarios. Consequently, distinct bidding curves are submitted to the intra-day market based on realized scenarios from the day-ahead market, as depicted in Fig. \ref{fig:Bid_IM}. 
Therefore, unlike in the day-ahead market, where bidding is primarily time-dependent, bidding in the intra-day market is influenced by the bidding curves submitted in the day-ahead market for each hour. This is because the bidding curve for the intra-day market is a second stage variable, determined after the clearance of the day-ahead market.
Notably, varying bidding curves for different scenarios underscore the sensitivity of second market bidding curves to fluctuations in first market prices. Moreover, the significance of appropriate clustering of electricity market prices becomes evident.

\begin{figure}[hbt!]
    \centering
    \subfloat[]{%
        \includegraphics[width=0.75\linewidth]{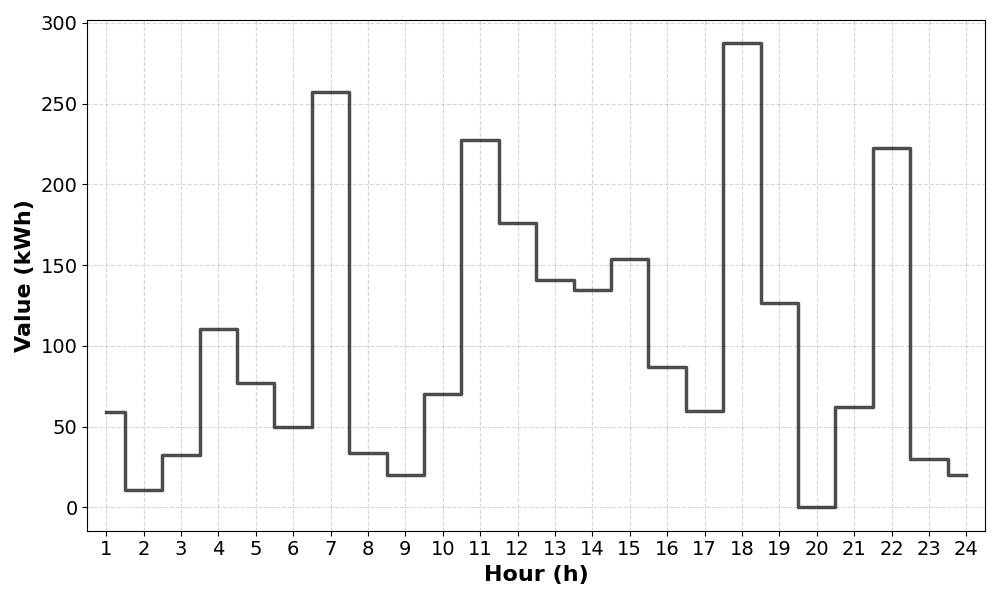}
        \label{fig:DAM}
    }
    \quad 
    \subfloat[]{%
        \includegraphics[width=0.75\linewidth]{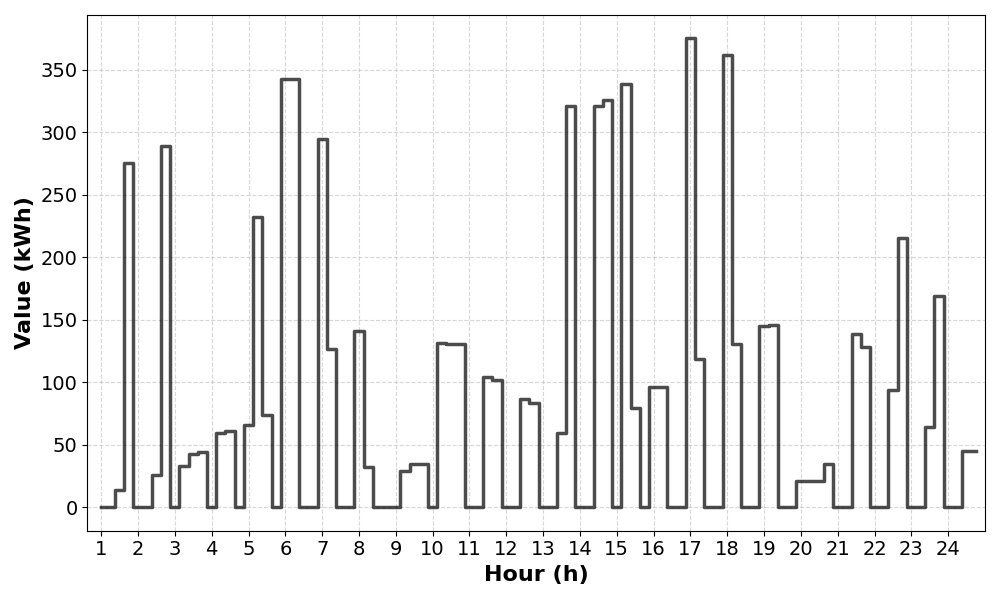}
        \label{fig:IDM}
    }
    \caption{Purchased power from electricity markets.  (a) Day-ahead market. (b) Intraday market.}
    \label{fig:market}
\end{figure}

The electricity purchased from the day-ahead and intraday markets, based on the average hourly day-ahead and quarter-hourly intraday prices, is illustrated in Fig. \ref{fig:market}. As depicted, the station buys power from the market offering the lowest prices during each period.

\begin{figure}[hbt!]
    \centering
    \subfloat[]{%
        \includegraphics[width=0.75\linewidth]{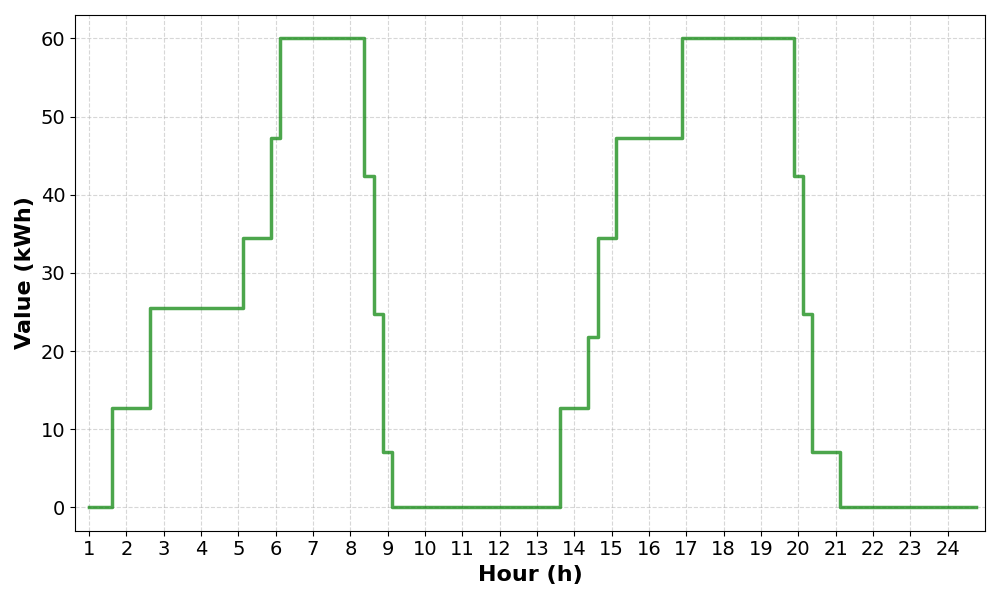}
        \label{fig:bl}
    }
    \quad 
    \subfloat[]{%
        \includegraphics[width=0.75\linewidth]{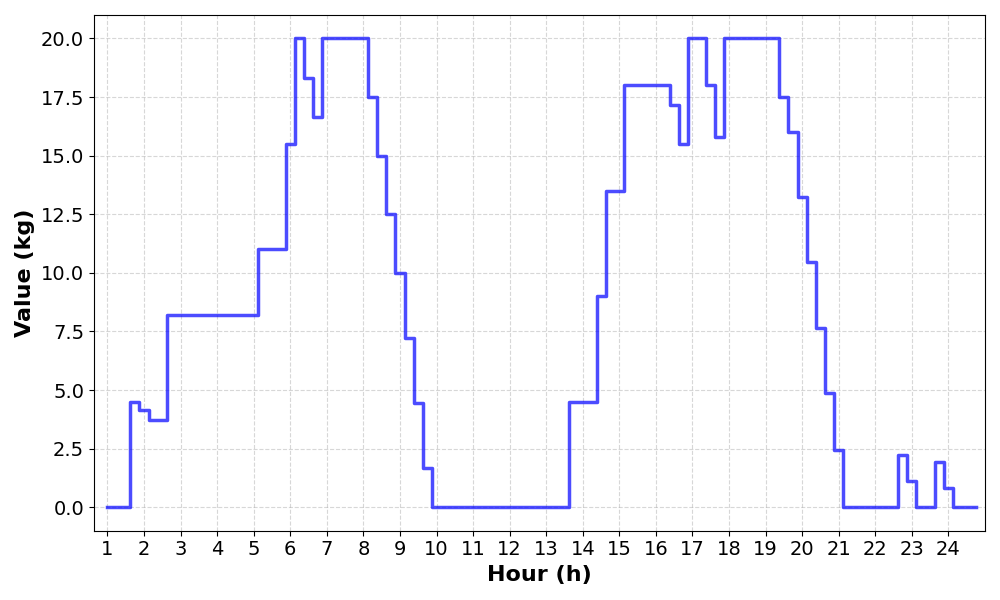}
        \label{fig:hl}
    }
    \caption{Stored amount. (a) Battery. (b) Tank.}
    \label{fig:level}
\end{figure}

The stored amounts of electricity and hydrogen in the battery and tank are shown in Fig. \ref{fig:level}. This reveals a pattern where the battery and tank are charged during off-peak hours and discharged during peak hours.

\begin{figure}[hbt!]
    \centering
\includegraphics[width=0.99\columnwidth]{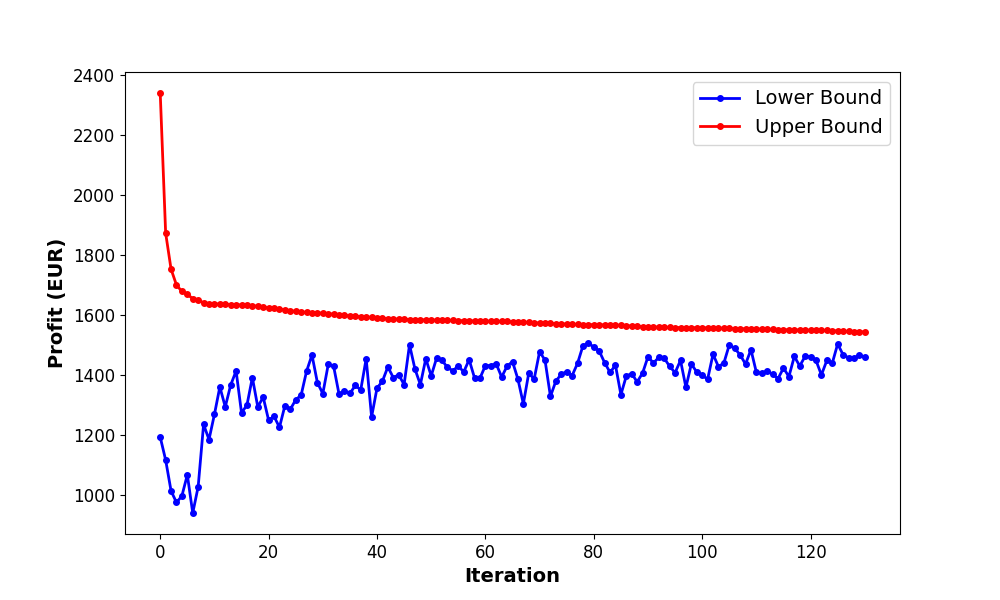}
    \caption{Convergence of our algorithm.}
    \label{fig:convergence}
\end{figure}

We address the challenges posed by the large scale and complexity of our study, which encompasses approximately 2,300 scenarios spanning 96 discrete 15-minute time intervals. This extensive scope exceeds the capabilities of traditional solvers, including CPLEX, Gurobi, and classic Benders, underscoring the necessity of our approach.
Fig.\ref{fig:convergence} demonstrates the convergence of our algorithm through numerous iterations. This emphasizes the efficiency of our method in tackling significant stochastic optimization challenges.

\begin{figure}[ht!]
    \centering
    \begin{subfigure}[b]{0.49\textwidth}
        \centering
        \includegraphics[width=\textwidth]{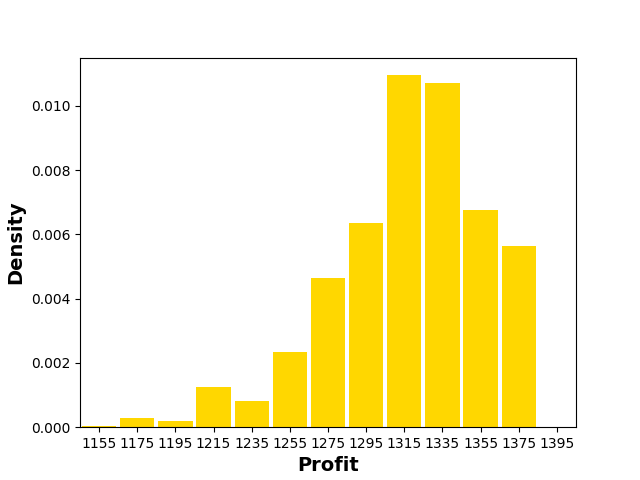}
        \caption{50\%}
        \label{fig:sub1}
    \end{subfigure}
    \hfill
    \begin{subfigure}[b]{0.49\textwidth}
        \centering
        \includegraphics[width=\textwidth]{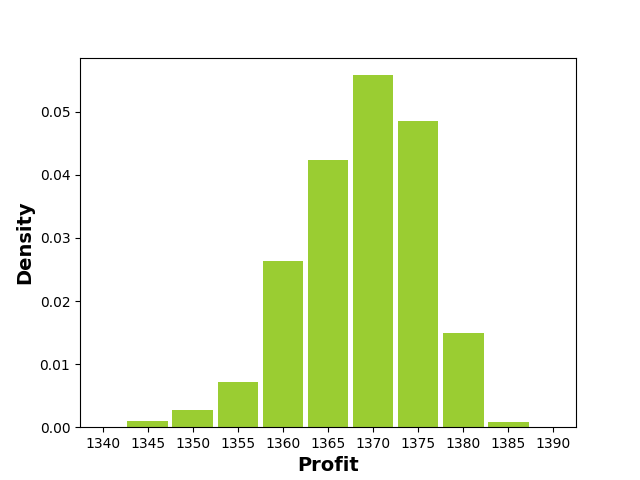}
        \caption{60\%}
        \label{fig:sub2}
    \end{subfigure}
    
    \begin{subfigure}[b]{0.49\textwidth}
        \centering
        \includegraphics[width=\textwidth]{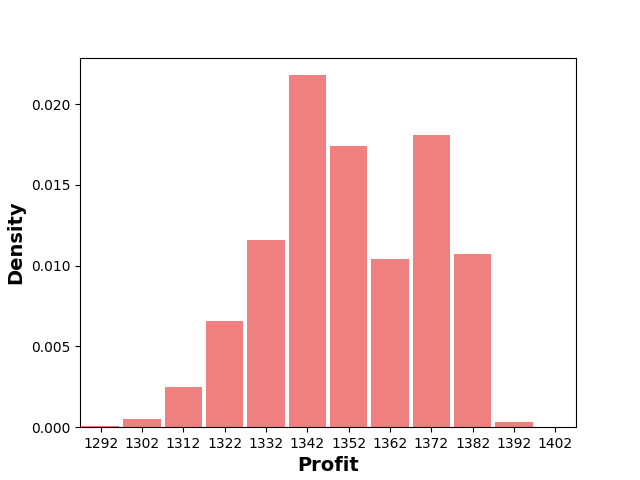}
        \caption{70\%}
        \label{fig:sub3}
    \end{subfigure}
    \hfill
    \begin{subfigure}[b]{0.49\textwidth}
        \centering
        \includegraphics[width=\textwidth]{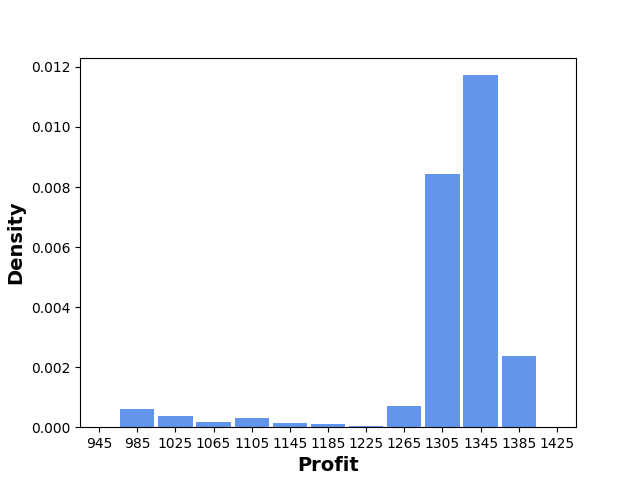}
        \caption{80\%}
        \label{fig:sub4}
    \end{subfigure}
    \caption{Histograms depicting the distribution of data at different confidence intervals.}
    \label{fig:subplots}
\end{figure}

\begin{figure}[hbt!]
    \centering
    \subfloat[]{%
        \includegraphics[width=0.55\linewidth]{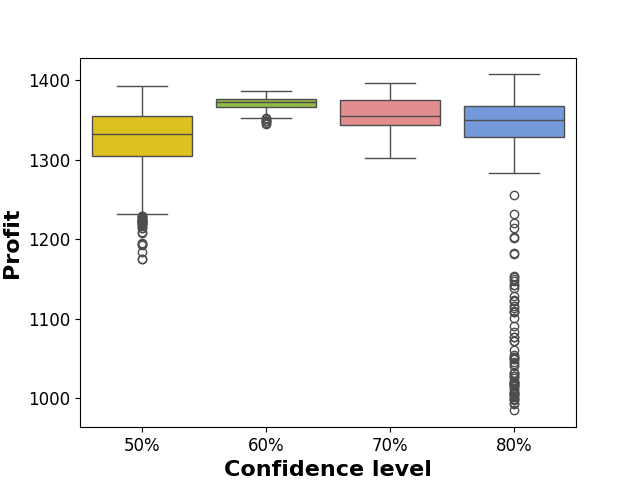}
        \label{fig:box}
    }
    \quad 
    \subfloat[]{%
        \includegraphics[width=0.55\linewidth]{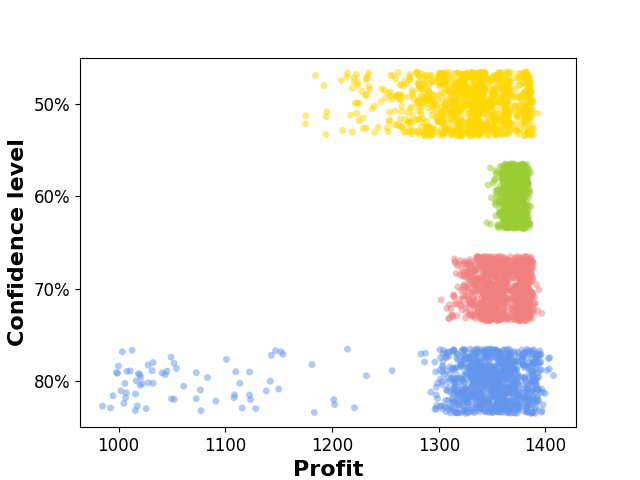}
        \label{fig:interval_value}
    }
    \caption{Analysis of various confidence intervals.}
    \label{fig:boxfish}
\end{figure}

Fig.\ref{fig:subplots} and Fig. \ref{fig:boxfish} depict the sensitivity analysis of the proposed model, aimed at determining the most advantageous range of prices for generating the discrete scenario set $\Xi$. Our primary objective is to ascertain which confidence levels of the price data obtained from the market yield the best results. This experiment is necessitated by the inherent complexity of scenario-based models, which often imposes limitations on the number of generated scenarios. Consequently, understanding the domain of the stochastic parameters most relevant for generating this limited number of scenarios becomes crucial. This exploration is imperative for the robustness and reliability of the proposed model under varying confidence levels. To address this, we investigate various confidence levels—namely 80\%, 70\%, 60\%, and 50\%, each representing a specific domain of electricity prices with $(1-\epsilon)\% $  frequency of most observed data (where $\epsilon$ equals 20\%, 30\%, 40\%, and 50\%). By simulating the model's responsiveness to different confidence levels, decision-makers can make well-informed decisions based on the model's output. In this test, three scenarios per hour are generated using the K-means method described in the section~\ref{K-means method}.

Fig.\ref{fig:subplots} and table \ref{tab:percentile} provide statistical metrics regarding the profits obtained from a Monte Carlo simulation comprising 1000 sets of randomly generated prices. These random prices are derived from the price distribution functions estimated using real market data. Although both the scenario set and the randomly generated prices in the Monte Carlo simulation are derived from the real market data, their specific values do not always match. This discrepancy arises because the scenario set provides a discrete representation of the continuous real market distribution, and actual realized prices may differ from the predefined scenarios. In such cases, the objective function in the Monte Carlo simulation is calculated using the optimal solution associated with the scenario that is closest to the selected price point. However, because this solution is optimized for the nearest scenario and not for the specific randomly selected price, it typically results in a lower profit compared to the profit associated with the scenario-based model.
In this simulation, four confidence intervals ranging from 50\% to 80\% are examined. 
Notably, at the 60\% confidence level, the dataset exhibits a mean value of 1371, representing the highest average profit among the intervals considered. Moreover, the variance, a measure of data dispersion, is notably low at this confidence level, registering at only 46. 
These findings suggest that in the 60\% confidence interval, the numbers are closer to the average, which means they're more consistent and precise. So, it seems like choosing the 60\% confidence interval is the best option for the scenario generation (considering three scenarios per hour) because it shows better consistency and less variation compared to others.

\begin{table}[ht!]
    \centering
    \begin{tabular}{ccccc} 
        \hline \hline
        & 50\% & 60\% & 70\% & 80\% \\ 
        \hline
        Average & 1327 & 1371 & 1357 & 1332 \\
        Variance & 1543 & 46 & 360 & 5914 \\
        \hline \hline
    \end{tabular}
    \caption{Average and Variance for Different Percentiles}
    \label{tab:percentile}
\end{table}
Regarding the results of our simulation, the box plots and individual value plots are employed for the statistical comparison of the proposed confidence levels. Fig.~\ref{fig:box} depicts a box plot that visualizes the median, quartiles, whisker, and potential outliers of the four confidence levels. Within the 60\% confidence level, the median is higher and the box length is narrower compared to other confidence levels. This indicates that data points cluster tightly around the median, implying greater consistency and precision. A shorter whisker length for the 60\% confidence level indicates less spread of the data outside the interquartile range (IQR), suggesting lower variability and higher precision. Comparing the outliers, which are data points that fall outside the whiskers of the box plot, shows that confidence levels 60\% and 70\% have fewer outliers, indicating less sensitivity of these confidence levels to price variabilities. Fig.~\ref{fig:interval_value} shows an individual value plot that visualizes the distribution spread of the profits obtained by the different confidence levels. This figure illustrates that the confidence level 60\% achieved a more favorable distribution with a lower variation. 

\section{Conclusion} \label{Conclusion}
In this paper, we propose a two-stage stochastic optimization model for a hybrid charging station, aiming to maximize owner profit while accommodating electric and hydrogen-powered vehicles. Energy procurement occurs through day-ahead and intraday auction-based electricity markets, where prices are uncertain. We utilize the k-means algorithm for stochastic scenario generation. Commercial solvers struggle with large-scale stochastic optimization problems due to their scope and stochastic nature. To address this, we employ a Benders decomposition-based algorithm within a two-stage stochastic programming framework. The Benders master problem determines day-ahead bidding decisions, while the sub-problem, utilizing SDDP, manages intraday bidding decisions and charging station scheduling.
We present numerical findings to demonstrate our contributions. Our results indicate that the charging station must implement diverse strategies to manage uncertainties effectively. Analysis reveals that revenue is most sensitive to changes in the selling price of hydrogen and least sensitive to electricity demand. Additionally, the multi-product charging station can effectively drive the advancement of electric and hydrogen vehicles, contributing to a more sustainable society.

\section*{Acknowledgments}
The research has been supported by the Czech Technical University in Prague, grant number SGS24/094/OHK5/2T/13, and the European Union under the project ROBOPROX (reg. no. CZ.02.01.01/00/22\_008/0004590). We also acknowledge the initial support provided by Mr. Anil Kaya.

\appendix

\section{}
\label{sec:appendix}

\begin{lemma} \label{lem:sp}
If we eliminate the binary variables from our proposed problem, the relaxed LP and the initial MILP models both yield identical optimal solutions.
\end{lemma}

\begin{proof}
Let's assume we remove the binary variables, which transforms the problem into a LP. This modification leads to simultaneous charging and discharging in both the battery and tank ($b_{c}\left(t,\xi\right), b_{d}\left(t,\xi\right) > 0$ and $h_{c}\left(t,\xi\right), h_{d}\left(t,\xi\right) > 0$). By applying the KKT optimality conditions, we adjust equations from \eqref{eq.2} to \eqref{eq.l2}, incorporating stationary conditions outlined in equations \eqref{eq.s1} through \eqref{eq.s4}.
\begin{flalign}
    \label{eq.l2}
    & m_{d}\left(t,\xi_1 \right) +m_{i}\left(t,\xi \right) + b_{d}\left(t,\xi\right) - l_{e}\left(t \right) - b_{c}\left(t,\xi\right) -  \mathrm{H} / \eta_{e} (h_{c}\left(t,\xi\right) -  h_{d}\left(t,\xi\right) + \notag\\
    &l_{h}\left(t \right))   = 0 \quad :\gamma\left(t,\xi\right) \quad \forall t, \forall \xi, \forall \xi_{1} \in \xi &&
    \end{flalign}
\begin{flalign}
    \label{eq.s1}
    & -\gamma\left(t,\xi\right) + \eta_{b} \tau\left(t,\xi\right) + \tau_{c,min}\left(t,\xi\right) - \tau_{c,max}\left(t,\xi\right) = 0  \quad \forall t, \forall \xi  &&
    \end{flalign}
\begin{flalign}
    \label{eq.s2}
    &  \gamma\left(t,\xi\right) - \tau\left(t,\xi\right) / \eta_{b} + \tau_{d,min}\left(t,\xi\right) - \tau_{d,max}\left(t,\xi\right) = 0 \quad \forall t, \forall \xi && 
    \end{flalign}
\begin{flalign}
    \label{eq.s3}
    & - \gamma\left(t,\xi\right) \mathrm{H}/ \eta_{e} + \eta_{h} \nu\left(t,\xi\right) + \nu_{c,min}\left(t,\xi\right) - \nu_{c,max}\left(t,\xi\right)= 0 \quad \forall t, \forall \xi && 
    \end{flalign}
\begin{flalign}
    \label{eq.s4}
    & \gamma\left(t,\xi\right) \mathrm{H}/ \eta_{e} - \nu\left(t,\xi\right) / \eta_{h} + \nu_{d,min}\left(t,\xi\right) - \nu_{d,max}\left(t,\xi\right)= 0 \quad \forall t, \forall \xi && 
    \end{flalign}

As the minimum charge and discharge capacities are zero, the Lagrangian multipliers $\tau_{c,min}\left(t,\xi\right)$, $\tau_{d,min}\left(t,\xi\right)$, $\nu_{c,min}\left(t,\xi\right)$, and $\nu_{d,min}\left(t,\xi\right)$ are all found to be zero. As a result, equations \eqref{eq.s5} and \eqref{eq.s6} emerge, ultimately leading us to deduce \eqref{eq.s7}.
\begin{flalign}
    \label{eq.s5}
    & (\gamma\left(t,\xi\right) + \tau_{c,max}\left(t,\xi\right) ) / \eta_{b} = \eta_{b} (\gamma\left(t,\xi\right) - \tau_{d,max}\left(t,\xi\right)) \quad \forall t, \forall \xi   && 
    \end{flalign}
\begin{flalign}
    \label{eq.s6}
    & ( \gamma\left(t,\xi\right) \mathrm{H}/ \eta_{e} + \nu_{c,max}\left(t,\xi\right) ) / \eta_{h} = \eta_{h} (\gamma\left(t,\xi\right) \mathrm{H}/ \eta_{e} - \nu_{d,max}\left(t,\xi\right)) \quad \forall t, \forall \xi   && 
    \end{flalign}
\begin{flalign}
    \label{eq.s7}
    & \gamma\left(t,\xi\right) 
    (\eta_{b} - 1/ \eta_{b}) + \gamma\left(t,\xi\right) \mathrm{H}/ \eta_{e} 
    (\eta_{h} - 1/ \eta_{h}) =  \tau_{c,max}\left(t,\xi\right) / \eta_{b} + \eta_{b} \tau_{d,max}\left(t,\xi\right) + \notag\\
    &\nu_{c,max}\left(t,\xi\right) / \eta_{h} +  \eta_{h} \nu_{d,max}\left(t,\xi\right) \quad \forall t, \forall \xi  &&
\end{flalign}

If we assume that the values of $b_{c}\left(t,\xi\right)$, $b_{d}\left(t,\xi\right)$, $h_{c}\left(t,\xi\right)$, and $h_{d}\left(t,\xi\right)$ are all greater than zero, then when we look at the right-hand side of equation \eqref{eq.s7}, it can only be either zero or strictly positive. However, the expression on the left-hand side is negative. This contradiction leads us to conclude that at least one of the variables $b_{c}\left(t,\xi\right)$ or $b_{d}\left(t,\xi\right)$ and $h_{c}\left(t,\xi\right)$ or $h_{d}\left(t,\xi\right)$ must be zero.
\end{proof}

\bibliographystyle{elsarticle-num-names} 

\bibliography{reference}

\begin{thebibliography}{44}
\expandafter\ifx\csname natexlab\endcsname\relax\def\natexlab#1{#1}\fi
\providecommand{\url}[1]{\texttt{#1}}
\providecommand{\href}[2]{#2}
\providecommand{\path}[1]{#1}
\providecommand{\DOIprefix}{doi:}
\providecommand{\ArXivprefix}{arXiv:}
\providecommand{\URLprefix}{URL: }
\providecommand{\Pubmedprefix}{pmid:}
\providecommand{\doi}[1]{\href{http://dx.doi.org/#1}{\path{#1}}}
\providecommand{\Pubmed}[1]{\href{pmid:#1}{\path{#1}}}
\providecommand{\bibinfo}[2]{#2}
\ifx\xfnm\relax \def\xfnm[#1]{\unskip,\space#1}\fi
\bibitem[{Dastjerdi et~al.(2023)Dastjerdi, Mosammam, Ahmadi, and Houshfar}]{dastjerdi2023transient}
\bibinfo{author}{S.~M. Dastjerdi}, \bibinfo{author}{Z.~M. Mosammam}, \bibinfo{author}{P.~Ahmadi}, \bibinfo{author}{E.~Houshfar},
\newblock \bibinfo{title}{Transient analysis and optimization of an off-grid hydrogen and electric vehicle charging station with temporary residences},
\newblock \bibinfo{journal}{Sustainable Cities and Society}  (\bibinfo{year}{2023}) \bibinfo{pages}{104742}.
\bibitem[{Elmasry et~al.(2024)Elmasry, Mansir, Abubakar, Ali, Aliyu, and Almamun}]{elmasry2024electricity}
\bibinfo{author}{Y.~Elmasry}, \bibinfo{author}{I.~B. Mansir}, \bibinfo{author}{Z.~Abubakar}, \bibinfo{author}{A.~Ali}, \bibinfo{author}{S.~Aliyu}, \bibinfo{author}{K.~Almamun},
\newblock \bibinfo{title}{Electricity-hydrogen nexus integrated with multi-level hydrogen storage, solar pv site, and electric-fuelcell car charging stations},
\newblock \bibinfo{journal}{International Journal of Hydrogen Energy}  (\bibinfo{year}{2024}).
\bibitem[{Sohrabi et~al.(2021)Sohrabi, Vahid-Pakdel, Mohammadi-Ivatloo, and Anvari-Moghaddam}]{sohrabi2021strategic}
\bibinfo{author}{F.~Sohrabi}, \bibinfo{author}{M.~Vahid-Pakdel}, \bibinfo{author}{B.~Mohammadi-Ivatloo}, \bibinfo{author}{A.~Anvari-Moghaddam},
\newblock \bibinfo{title}{Strategic planning of power to gas energy storage facilities in electricity market},
\newblock \bibinfo{journal}{Sustainable Energy Technologies and Assessments} \bibinfo{volume}{46} (\bibinfo{year}{2021}) \bibinfo{pages}{101238}.
\bibitem[{Cai et~al.(2023)Cai, Mi, Ma, Li, Li, and Wang}]{cai2023hierarchical}
\bibinfo{author}{P.~Cai}, \bibinfo{author}{Y.~Mi}, \bibinfo{author}{S.~Ma}, \bibinfo{author}{H.~Li}, \bibinfo{author}{D.~Li}, \bibinfo{author}{P.~Wang},
\newblock \bibinfo{title}{Hierarchical game for integrated energy system and electricity-hydrogen hybrid charging station under distributionally robust optimization},
\newblock \bibinfo{journal}{Energy} \bibinfo{volume}{283} (\bibinfo{year}{2023}) \bibinfo{pages}{128471}.
\bibitem[{Fang et~al.(2023)Fang, Wang, Dong, Yang, and Sun}]{fang2023optimal}
\bibinfo{author}{X.~Fang}, \bibinfo{author}{Y.~Wang}, \bibinfo{author}{W.~Dong}, \bibinfo{author}{Q.~Yang}, \bibinfo{author}{S.~Sun},
\newblock \bibinfo{title}{Optimal energy management of multiple electricity-hydrogen integrated charging stations},
\newblock \bibinfo{journal}{Energy} \bibinfo{volume}{262} (\bibinfo{year}{2023}) \bibinfo{pages}{125624}.
\bibitem[{S{\'a}nchez-S{\'a}inz et~al.(2019)S{\'a}nchez-S{\'a}inz, Garc{\'\i}a-V{\'a}zquez, Llorens~Iborra, and Fern{\'a}ndez-Ram{\'\i}rez}]{sanchez2019methodology}
\bibinfo{author}{H.~S{\'a}nchez-S{\'a}inz}, \bibinfo{author}{C.-A. Garc{\'\i}a-V{\'a}zquez}, \bibinfo{author}{F.~Llorens~Iborra}, \bibinfo{author}{L.~M. Fern{\'a}ndez-Ram{\'\i}rez},
\newblock \bibinfo{title}{Methodology for the optimal design of a hybrid charging station of electric and fuel cell vehicles supplied by renewable energies and an energy storage system},
\newblock \bibinfo{journal}{Sustainability} \bibinfo{volume}{11} (\bibinfo{year}{2019}) \bibinfo{pages}{5743}.
\bibitem[{Raj et~al.(2022)Raj, Saravanan, Durairaj, Selvakumar, Sagar, Ganavi, and Pavithra}]{raj2022energy}
\bibinfo{author}{E.~Raj}, \bibinfo{author}{A.~Saravanan}, \bibinfo{author}{U.~Durairaj}, \bibinfo{author}{S.~Selvakumar}, \bibinfo{author}{B.~Sagar}, \bibinfo{author}{M.~Ganavi}, \bibinfo{author}{G.~Pavithra},
\newblock \bibinfo{title}{Energy management system for charging stations for electric and hydrogen vehicles using solar, hydrogen and fuel cell technology},
\newblock in: \bibinfo{booktitle}{AIP Conference Proceedings}, volume \bibinfo{volume}{2519}, \bibinfo{organization}{AIP Publishing}, \bibinfo{year}{2022}.
\bibitem[{Sriyakul and Jermsittiparsert(2021)}]{sriyakul2021risk}
\bibinfo{author}{T.~Sriyakul}, \bibinfo{author}{K.~Jermsittiparsert},
\newblock \bibinfo{title}{Risk-constrained design of autonomous hybrid refueling station for hydrogen and electric vehicles using information gap decision theory},
\newblock \bibinfo{journal}{International Journal of Hydrogen Energy} \bibinfo{volume}{46} (\bibinfo{year}{2021}) \bibinfo{pages}{1682--1693}.
\bibitem[{Sohrabi et~al.(2022)Sohrabi, Rohaninejad, Hesamzadeh, and Bem{\v{s}}}]{sohrabi2022optimal}
\bibinfo{author}{F.~Sohrabi}, \bibinfo{author}{M.~Rohaninejad}, \bibinfo{author}{M.~R. Hesamzadeh}, \bibinfo{author}{J.~Bem{\v{s}}},
\newblock \bibinfo{title}{Optimal trading of a hybrid electric, hydrogen and gas fueling station in day-ahead and intra-day markets: Modeling aspect},
\newblock in: \bibinfo{booktitle}{International Conference on Operations Research}, \bibinfo{organization}{Springer}, \bibinfo{year}{2022}, pp. \bibinfo{pages}{289--295}.
\bibitem[{Guo and Gong(2022)}]{guo2022energy}
\bibinfo{author}{G.~Guo}, \bibinfo{author}{Y.~Gong},
\newblock \bibinfo{title}{Energy management of intelligent solar parking lot with ev charging and fcev refueling based on deep reinforcement learning},
\newblock \bibinfo{journal}{International Journal of Electrical Power \& Energy Systems} \bibinfo{volume}{140} (\bibinfo{year}{2022}) \bibinfo{pages}{108061}.
\bibitem[{{\c{C}}i{\c{c}}ek(2022)}]{cciccek2022optimal}
\bibinfo{author}{A.~{\c{C}}i{\c{c}}ek},
\newblock \bibinfo{title}{Optimal operation of an all-in-one ev station with photovoltaic system including charging, battery swapping and hydrogen refueling},
\newblock \bibinfo{journal}{International Journal of Hydrogen Energy} \bibinfo{volume}{47} (\bibinfo{year}{2022}) \bibinfo{pages}{32405--32424}.
\bibitem[{Xu et~al.(2022)Xu, Hu, Liu, Du, Huang, and Chen}]{xu2022robust}
\bibinfo{author}{X.~Xu}, \bibinfo{author}{W.~Hu}, \bibinfo{author}{W.~Liu}, \bibinfo{author}{Y.~Du}, \bibinfo{author}{Q.~Huang}, \bibinfo{author}{Z.~Chen},
\newblock \bibinfo{title}{Robust energy management for an on-grid hybrid hydrogen refueling and battery swapping station based on renewable energy},
\newblock \bibinfo{journal}{Journal of Cleaner Production} \bibinfo{volume}{331} (\bibinfo{year}{2022}) \bibinfo{pages}{129954}.
\bibitem[{Schr{\"o}der et~al.(2020)Schr{\"o}der, Abdin, and M{\'e}rida}]{schroder2020optimization}
\bibinfo{author}{M.~Schr{\"o}der}, \bibinfo{author}{Z.~Abdin}, \bibinfo{author}{W.~M{\'e}rida},
\newblock \bibinfo{title}{Optimization of distributed energy resources for electric vehicle charging and fuel cell vehicle refueling},
\newblock \bibinfo{journal}{Applied energy} \bibinfo{volume}{277} (\bibinfo{year}{2020}) \bibinfo{pages}{115562}.
\bibitem[{Mehrjerdi(2019)}]{mehrjerdi2019off}
\bibinfo{author}{H.~Mehrjerdi},
\newblock \bibinfo{title}{Off-grid solar powered charging station for electric and hydrogen vehicles including fuel cell and hydrogen storage},
\newblock \bibinfo{journal}{International journal of hydrogen Energy} \bibinfo{volume}{44} (\bibinfo{year}{2019}) \bibinfo{pages}{11574--11583}.
\bibitem[{Sohrabi et~al.(2023)Sohrabi, Rohaninejad, Hesamzadeh, and Bem{\v{s}}}]{sohrabi2023coordinated}
\bibinfo{author}{F.~Sohrabi}, \bibinfo{author}{M.~Rohaninejad}, \bibinfo{author}{M.~R. Hesamzadeh}, \bibinfo{author}{J.~Bem{\v{s}}},
\newblock \bibinfo{title}{Coordinated bidding of multi-product charging station in electricity markets using rolling planning and sample average approximation},
\newblock \bibinfo{journal}{International Journal of Electrical Power \& Energy Systems} \bibinfo{volume}{146} (\bibinfo{year}{2023}) \bibinfo{pages}{108786}.
\bibitem[{Wang et~al.(2020)Wang, Kazemi, Nojavan, and Jermsittiparsert}]{wang2020robust}
\bibinfo{author}{Y.~Wang}, \bibinfo{author}{M.~Kazemi}, \bibinfo{author}{S.~Nojavan}, \bibinfo{author}{K.~Jermsittiparsert},
\newblock \bibinfo{title}{Robust design of off-grid solar-powered charging station for hydrogen and electric vehicles via robust optimization approach},
\newblock \bibinfo{journal}{International Journal of Hydrogen Energy} \bibinfo{volume}{45} (\bibinfo{year}{2020}) \bibinfo{pages}{18995--19006}.
\bibitem[{Ampah et~al.(2022)Ampah, Afrane, Agyekum, Adun, Yusuf, and Bamisile}]{ampah2022electric}
\bibinfo{author}{J.~D. Ampah}, \bibinfo{author}{S.~Afrane}, \bibinfo{author}{E.~B. Agyekum}, \bibinfo{author}{H.~Adun}, \bibinfo{author}{A.~A. Yusuf}, \bibinfo{author}{O.~Bamisile},
\newblock \bibinfo{title}{Electric vehicles development in sub-saharan africa: Performance assessment of standalone renewable energy systems for hydrogen refuelling and electricity charging stations (hrecs)},
\newblock \bibinfo{journal}{Journal of Cleaner Production} \bibinfo{volume}{376} (\bibinfo{year}{2022}) \bibinfo{pages}{134238}.
\bibitem[{Massana et~al.(2022)Massana, Burgas, Herraiz, Colomer, and Pous}]{massana2022multi}
\bibinfo{author}{J.~Massana}, \bibinfo{author}{L.~Burgas}, \bibinfo{author}{S.~Herraiz}, \bibinfo{author}{J.~Colomer}, \bibinfo{author}{C.~Pous},
\newblock \bibinfo{title}{Multi-vector energy management system including scheduling electrolyser, electric vehicle charging station and other assets in a real scenario},
\newblock \bibinfo{journal}{Journal of Cleaner Production} \bibinfo{volume}{380} (\bibinfo{year}{2022}) \bibinfo{pages}{134996}.
\bibitem[{Zeng et~al.(2023)Zeng, Wang, Zhang, Wang, Tang, and Wang}]{zeng2023optimal}
\bibinfo{author}{B.~Zeng}, \bibinfo{author}{W.~Wang}, \bibinfo{author}{W.~Zhang}, \bibinfo{author}{Y.~Wang}, \bibinfo{author}{C.~Tang}, \bibinfo{author}{J.~Wang},
\newblock \bibinfo{title}{Optimal configuration planning of vehicle sharing station-based electro-hydrogen micro-energy systems for transportation decarbonization},
\newblock \bibinfo{journal}{Journal of Cleaner Production} \bibinfo{volume}{387} (\bibinfo{year}{2023}) \bibinfo{pages}{135906}.
\bibitem[{Silva-Rodriguez et~al.(2022)Silva-Rodriguez, Sanjab, Fumagalli, Virag, and Gibescu}]{silva2022short}
\bibinfo{author}{L.~Silva-Rodriguez}, \bibinfo{author}{A.~Sanjab}, \bibinfo{author}{E.~Fumagalli}, \bibinfo{author}{A.~Virag}, \bibinfo{author}{M.~Gibescu},
\newblock \bibinfo{title}{Short term wholesale electricity market designs: A review of identified challenges and promising solutions},
\newblock \bibinfo{journal}{Renewable and Sustainable Energy Reviews} \bibinfo{volume}{160} (\bibinfo{year}{2022}) \bibinfo{pages}{112228}.
\bibitem[{Braun and Brunner(2018)}]{braun2018price}
\bibinfo{author}{S.~M. Braun}, \bibinfo{author}{C.~Brunner},
\newblock \bibinfo{title}{Price sensitivity of hourly day-ahead and quarter-hourly intraday auctions in germany},
\newblock \bibinfo{journal}{Zeitschrift f{\"u}r Energiewirtschaft} \bibinfo{volume}{42} (\bibinfo{year}{2018}) \bibinfo{pages}{257--270}.
\bibitem[{Ocker and Ehrhart(2017)}]{ocker2017german}
\bibinfo{author}{F.~Ocker}, \bibinfo{author}{K.-M. Ehrhart},
\newblock \bibinfo{title}{The “german paradox” in the balancing power markets},
\newblock \bibinfo{journal}{Renewable and Sustainable Energy Reviews} \bibinfo{volume}{67} (\bibinfo{year}{2017}) \bibinfo{pages}{892--898}.
\bibitem[{Finnah et~al.(2022)Finnah, G{\"o}nsch, and Ziel}]{finnah2022integrated}
\bibinfo{author}{B.~Finnah}, \bibinfo{author}{J.~G{\"o}nsch}, \bibinfo{author}{F.~Ziel},
\newblock \bibinfo{title}{Integrated day-ahead and intraday self-schedule bidding for energy storage systems using approximate dynamic programming},
\newblock \bibinfo{journal}{European Journal of Operational Research} \bibinfo{volume}{301} (\bibinfo{year}{2022}) \bibinfo{pages}{726--746}.
\bibitem[{MacQueen et~al.(1967)}]{macqueen1967some}
\bibinfo{author}{J.~MacQueen}, et~al.,
\newblock \bibinfo{title}{Some methods for classification and analysis of multivariate observations},
\newblock in: \bibinfo{booktitle}{Proceedings of the fifth Berkeley symposium on mathematical statistics and probability}, volume~\bibinfo{volume}{1}, \bibinfo{organization}{Oakland, CA, USA}, \bibinfo{year}{1967}, pp. \bibinfo{pages}{281--297}.
\bibitem[{Seljom et~al.(2021)Seljom, Kvalbein, Hellemo, Kaut, and Ortiz}]{seljom2021stochastic}
\bibinfo{author}{P.~Seljom}, \bibinfo{author}{L.~Kvalbein}, \bibinfo{author}{L.~Hellemo}, \bibinfo{author}{M.~Kaut}, \bibinfo{author}{M.~M. Ortiz},
\newblock \bibinfo{title}{Stochastic modelling of variable renewables in long-term energy models: Dataset, scenario generation \& quality of results},
\newblock \bibinfo{journal}{Energy} \bibinfo{volume}{236} (\bibinfo{year}{2021}) \bibinfo{pages}{121415}.
\bibitem[{Lloyd(1982)}]{lloyd1982least}
\bibinfo{author}{S.~Lloyd},
\newblock \bibinfo{title}{Least squares quantization in pcm},
\newblock \bibinfo{journal}{IEEE transactions on information theory} \bibinfo{volume}{28} (\bibinfo{year}{1982}) \bibinfo{pages}{129--137}.
\bibitem[{Yaghoubi-Nia et~al.(2021)Yaghoubi-Nia, Hashemi-Dezaki, and Niasar}]{yaghoubi2021optimal}
\bibinfo{author}{M.-R. Yaghoubi-Nia}, \bibinfo{author}{H.~Hashemi-Dezaki}, \bibinfo{author}{A.~H. Niasar},
\newblock \bibinfo{title}{Optimal stochastic scenario-based allocation of smart grids’ renewable and non-renewable distributed generation units and protective devices},
\newblock \bibinfo{journal}{Sustainable Energy Technologies and Assessments} \bibinfo{volume}{44} (\bibinfo{year}{2021}) \bibinfo{pages}{101033}.
\bibitem[{Rebennack(2014)}]{rebennack2014generation}
\bibinfo{author}{S.~Rebennack},
\newblock \bibinfo{title}{Generation expansion planning under uncertainty with emissions quotas},
\newblock \bibinfo{journal}{Electric Power Systems Research} \bibinfo{volume}{114} (\bibinfo{year}{2014}) \bibinfo{pages}{78--85}.
\bibitem[{BnnoBRs(1962)}]{bnnobrs1962partitioning}
\bibinfo{author}{J.~BnnoBRs},
\newblock \bibinfo{title}{Partitioning procedures for solving mixed-variables programming problems},
\newblock \bibinfo{journal}{Numer. Math} \bibinfo{volume}{4} (\bibinfo{year}{1962}) \bibinfo{pages}{238--252}.
\bibitem[{Rahmaniani et~al.(2017)Rahmaniani, Crainic, Gendreau, and Rei}]{rahmaniani2017benders}
\bibinfo{author}{R.~Rahmaniani}, \bibinfo{author}{T.~G. Crainic}, \bibinfo{author}{M.~Gendreau}, \bibinfo{author}{W.~Rei},
\newblock \bibinfo{title}{The benders decomposition algorithm: A literature review},
\newblock \bibinfo{journal}{European Journal of Operational Research} \bibinfo{volume}{259} (\bibinfo{year}{2017}) \bibinfo{pages}{801--817}.
\bibitem[{Alizadeh and Jadid(2015)}]{alizadeh2015dynamic}
\bibinfo{author}{B.~Alizadeh}, \bibinfo{author}{S.~Jadid},
\newblock \bibinfo{title}{A dynamic model for coordination of generation and transmission expansion planning in power systems},
\newblock \bibinfo{journal}{International Journal of Electrical Power \& Energy Systems} \bibinfo{volume}{65} (\bibinfo{year}{2015}) \bibinfo{pages}{408--418}.
\bibitem[{Rahmani and Mahoodian(2017)}]{rahmani2017strategic}
\bibinfo{author}{D.~Rahmani}, \bibinfo{author}{V.~Mahoodian},
\newblock \bibinfo{title}{Strategic and operational supply chain network design to reduce carbon emission considering reliability and robustness},
\newblock \bibinfo{journal}{Journal of Cleaner Production} \bibinfo{volume}{149} (\bibinfo{year}{2017}) \bibinfo{pages}{607--620}.
\bibitem[{Pishvaee et~al.(2014)Pishvaee, Razmi, and Torabi}]{pishvaee2014accelerated}
\bibinfo{author}{M.~S. Pishvaee}, \bibinfo{author}{J.~Razmi}, \bibinfo{author}{S.~A. Torabi},
\newblock \bibinfo{title}{An accelerated benders decomposition algorithm for sustainable supply chain network design under uncertainty: A case study of medical needle and syringe supply chain},
\newblock \bibinfo{journal}{Transportation Research Part E: Logistics and Transportation Review} \bibinfo{volume}{67} (\bibinfo{year}{2014}) \bibinfo{pages}{14--38}.
\bibitem[{Pereira and Pinto(1991)}]{pereira1991multi}
\bibinfo{author}{M.~V. Pereira}, \bibinfo{author}{L.~M. Pinto},
\newblock \bibinfo{title}{Multi-stage stochastic optimization applied to energy planning},
\newblock \bibinfo{journal}{Mathematical programming} \bibinfo{volume}{52} (\bibinfo{year}{1991}) \bibinfo{pages}{359--375}.
\bibitem[{Lei et~al.(2024)Lei, Huang, Xu, Zhu, Yang, Liu, and Hu}]{lei2024optimal}
\bibinfo{author}{Q.~Lei}, \bibinfo{author}{Y.~Huang}, \bibinfo{author}{X.~Xu}, \bibinfo{author}{F.~Zhu}, \bibinfo{author}{Y.~Yang}, \bibinfo{author}{J.~Liu}, \bibinfo{author}{W.~Hu},
\newblock \bibinfo{title}{Optimal scheduling of a renewable energy-based park power system: a novel hybrid sddp/mpc approach},
\newblock \bibinfo{journal}{International Journal of Electrical Power \& Energy Systems} \bibinfo{volume}{157} (\bibinfo{year}{2024}) \bibinfo{pages}{109892}.
\bibitem[{Rebennack(2016)}]{rebennack2016combining}
\bibinfo{author}{S.~Rebennack},
\newblock \bibinfo{title}{Combining sampling-based and scenario-based nested benders decomposition methods: application to stochastic dual dynamic programming},
\newblock \bibinfo{journal}{Mathematical Programming} \bibinfo{volume}{156} (\bibinfo{year}{2016}) \bibinfo{pages}{343--389}.
\bibitem[{Ding et~al.(2021)Ding, Zhang, Lu, Qu, Shahidehpour, He, and Chen}]{ding2021multi}
\bibinfo{author}{T.~Ding}, \bibinfo{author}{X.~Zhang}, \bibinfo{author}{R.~Lu}, \bibinfo{author}{M.~Qu}, \bibinfo{author}{M.~Shahidehpour}, \bibinfo{author}{Y.~He}, \bibinfo{author}{T.~Chen},
\newblock \bibinfo{title}{Multi-stage distributionally robust stochastic dual dynamic programming to multi-period economic dispatch with virtual energy storage},
\newblock \bibinfo{journal}{IEEE Transactions on Sustainable Energy} \bibinfo{volume}{13} (\bibinfo{year}{2021}) \bibinfo{pages}{146--158}.
\bibitem[{Steeger and Rebennack(2017)}]{steeger2017dynamic}
\bibinfo{author}{G.~Steeger}, \bibinfo{author}{S.~Rebennack},
\newblock \bibinfo{title}{Dynamic convexification within nested benders decomposition using lagrangian relaxation: An application to the strategic bidding problem},
\newblock \bibinfo{journal}{European Journal of Operational Research} \bibinfo{volume}{257} (\bibinfo{year}{2017}) \bibinfo{pages}{669--686}.
\bibitem[{Fleten and Pettersen(2005)}]{fleten2005constructing}
\bibinfo{author}{S.-E. Fleten}, \bibinfo{author}{E.~Pettersen},
\newblock \bibinfo{title}{Constructing bidding curves for a price-taking retailer in the norwegian electricity market},
\newblock \bibinfo{journal}{IEEE Transactions on Power Systems} \bibinfo{volume}{20} (\bibinfo{year}{2005}) \bibinfo{pages}{701--708}.
\bibitem[{Zeng et~al.(2020)Zeng, Li, and Wang}]{zeng2020scenario}
\bibinfo{author}{Y.~Zeng}, \bibinfo{author}{C.~Li}, \bibinfo{author}{H.~Wang},
\newblock \bibinfo{title}{Scenario-set-based economic dispatch of power system with wind power and energy storage system},
\newblock \bibinfo{journal}{IEEE Access} \bibinfo{volume}{8} (\bibinfo{year}{2020}) \bibinfo{pages}{109105--109119}.
\bibitem[{Noorollahi et~al.(2022)Noorollahi, Golshanfard, and Hashemi-Dezaki}]{noorollahi2022scenario}
\bibinfo{author}{Y.~Noorollahi}, \bibinfo{author}{A.~Golshanfard}, \bibinfo{author}{H.~Hashemi-Dezaki},
\newblock \bibinfo{title}{A scenario-based approach for optimal operation of energy hub under different schemes and structures},
\newblock \bibinfo{journal}{Energy} \bibinfo{volume}{251} (\bibinfo{year}{2022}) \bibinfo{pages}{123740}.
\bibitem[{Ashraf et~al.(2017)Ashraf, Gupta, Choudhary, and Chakrabarti}]{ashraf2017voltage}
\bibinfo{author}{S.~M. Ashraf}, \bibinfo{author}{A.~Gupta}, \bibinfo{author}{D.~K. Choudhary}, \bibinfo{author}{S.~Chakrabarti},
\newblock \bibinfo{title}{Voltage stability monitoring of power systems using reduced network and artificial neural network},
\newblock \bibinfo{journal}{International Journal of Electrical Power \& Energy Systems} \bibinfo{volume}{87} (\bibinfo{year}{2017}) \bibinfo{pages}{43--51}.
\bibitem[{Fatouros et~al.(2017)Fatouros, Konstantelos, Papadaskalopoulos, and Strbac}]{fatouros2017stochastic}
\bibinfo{author}{P.~Fatouros}, \bibinfo{author}{I.~Konstantelos}, \bibinfo{author}{D.~Papadaskalopoulos}, \bibinfo{author}{G.~Strbac},
\newblock \bibinfo{title}{Stochastic dual dynamic programming for operation of der aggregators under multi-dimensional uncertainty},
\newblock \bibinfo{journal}{IEEE transactions on sustainable energy} \bibinfo{volume}{10} (\bibinfo{year}{2017}) \bibinfo{pages}{459--469}.
\bibitem[{{ENTSO-E}(2024)}]{entsoe}
\bibinfo{author}{{ENTSO-E}}, \bibinfo{title}{{European Network of Transmission System Operators for Electricity}}, \bibinfo{howpublished}{\url{https://www.entsoe.eu/}}, \bibinfo{year}{2024}.

\end{thebibliography}

\end{document}